\documentclass[12pt,a4paper,oneside]{amsart}

\usepackage{amsmath}
\usepackage{amsthm}
\usepackage{amsfonts}
\usepackage{amssymb}
\usepackage{amscd}

\usepackage[T1]{fontenc}
\usepackage{textcomp}

\usepackage[warnunknown]{ucs}
\usepackage[utf8x]{inputenc}

\usepackage{mathpazo}
\usepackage[scaled]{helvet}

\usepackage{calc}

\usepackage[pdftex,final]{graphicx}
\usepackage[usenames]{color}

\usepackage{float}

\usepackage{microtype}

\usepackage{url}

\usepackage[verbose=true,heightrounded=true]{geometry}

\numberwithin{equation}{section}

\begin{document}
\title[Plane puzzles]{Periodical plane puzzles with numbers}
\author{Jorge Rezende}
\address{Grupo de F\'{i}sica-Matem\'{a}tica da Universidade de Lisboa, Av. Prof. Gama
Pinto 2, 1649-003 Lisboa, Portugal, and Departamento de Matem\'{a}tica,
Faculdade de Ci\^{e}ncias da Universidade de Lisboa}
\email{rezende@cii.fc.ul.pt}
\thanks{The Mathematical Physics Group is supported by the (Portuguese) Foundation
for Science and Technology (FCT). }
\thanks{This paper is in final form and no version of it will be submitted for
publication elsewhere.}
\maketitle

\section{Introduction}

Consider a periodical (in two independent directions) tiling of the plane
with polygons (faces). In this article we shall only give examples using
squares, regular hexagons, equilateral triangles and parallelograms
(``unions'' of two equilateral triangles). We shall call some ``multiple''
of the fundamental region ``the board''. We naturally identify pairs of
corresponding edges of the the board. Figures 9 and 19--29 show different
boards. The ``border'' of the board is represented by a yellow thick line,
unless part of it or all of it is the edge of a face.

The board is tiled by a finite number of polygons. Construct polygonal
plates in the same number, shape and size as the polygons of the board.
Adjacent to each side of each plate draw a number, or two numbers, like it
is shown in Figures 1 and 18-29. Figure 1 shows the obvious possibility of
having plates with simple drawings, coloured drawings, etc.
\begin{figure}[h]
  \centering
  \includegraphics[width=3.8692in]{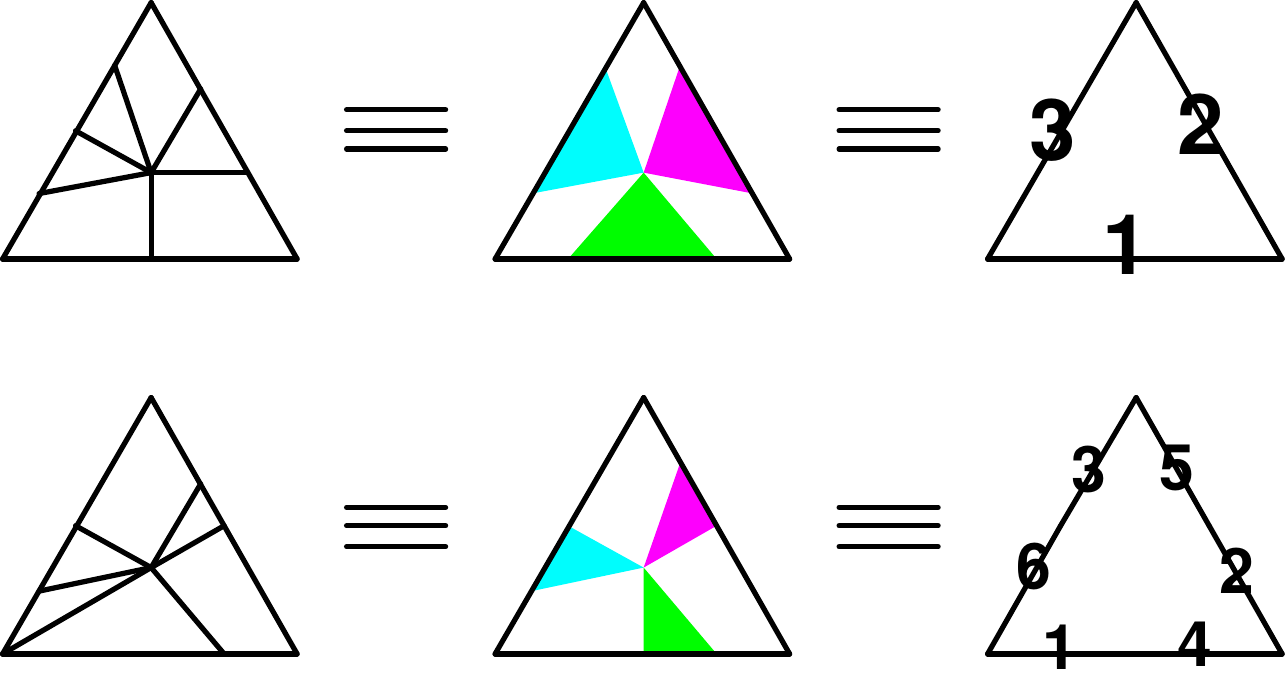}  
  \caption{}
  \label{fig:01}
\end{figure}

Now the game is to put the plates over the board polygons in such a way that
the numbers near each board edge are equal. If there is at least one
solution of this puzzle one says that we have a periodical plane puzzle with
numbers.

These puzzles are a tool in teaching and learning mathematics. For those
that already have some mathematical knowledge, they are a source for many
examples and exercises, that go from the elementary to complex ones, in
combinatorics, group theory (including symmetry and permutation groups),
programming, and so on. The object of this work is to point out some
possibilities by giving simple examples.

The computer is the only practical way of ``materializing'' infinite
periodical plane puzzles. Hence, these puzzles can be very well put into
practice as computer games.

This article follows some others on \textit{puzzles with numbers}.
See \cite{rezende1}, \cite{rezende2}, \cite{rezende3}.

\section{Plane symmetries}

A function $\omega :\mathbb{R}^2\rightarrow \mathbb{R}^2$ is an \textit{isometry}
if it preserves the Euclidean distance: $\left\| \omega \left( \alpha
\right) -\omega \left( \beta \right) \right\| =\left\| \alpha -\beta
\right\| $, for every $\alpha ,\beta \in \mathbb{R}^2$. We denote $\mathcal{I}$
the group of isometries of $\mathbb{R}^2$.

Every \textit{translation} $\omega :\alpha \longmapsto u+\alpha $, for some 
$u\in \mathbb{R}^2$, is an isometry.

Every linear isometry belongs to $O_2$, the \textit{orthogonal group}. We
denote the identity by $i$.

Every isometry is a composition of a translation with an orthogonal
transformation, i.e., it is of the form $\alpha \longmapsto u+\eta \left(
\alpha \right) $, for some $u\in \mathbb{R}^2$ and $\eta \in O_2$. The pair 
$\left( u,\eta \right) $ defines the isometry. We represent the isometry by 
$\omega \equiv \left( u,\eta \right) \in \mathcal{I}$.

A \textit{rotation} about the point $\gamma $ is an isometry $\omega $,
$\omega \left( \alpha \right) =u+\eta \left( \alpha \right) $, if $\det \eta
=1$ and $\omega \left( \gamma \right) =\gamma $, i.e., $u=\gamma -\eta
\left( \gamma \right) $. Hence, if $\eta \neq i$, $\gamma =\left( i-\eta
\right) ^{-1}\left( u\right) $. Notice that, if $\det \eta =1$ and $\eta
\neq i$, then $i-\eta $ is invertible.

If $\det \eta =1$, one says that we have a \textit{direct isometry}. Every
direct isometry is a translation or a rotation\cite{Armstrong}. The
identity, $\omega =i$, is a rotation and a translation. If $\omega \neq i$
and $\eta =i$, then $\omega $ is a translation. If $\omega \neq i$ and $\eta
\neq i$, then $\omega $ is a rotation.

If $\det \eta =-1$, every $\gamma \in \mathbb{R}^2$ can be decomposed $\gamma
=\gamma _1+\gamma _2$, with $\gamma _1\perp \gamma _2$, such that $\eta
\left( \gamma \right) =\eta \left( \gamma _1+\gamma _2\right) =-\gamma
_1+\gamma _2$. The isometry $\omega \left( \alpha \right) =u+\eta \left(
\alpha \right) $ is of the form $\omega \left( \alpha _1+\alpha _2\right)
=u_1+u_2-\alpha _1+\alpha _2$. If $u_2=0$, $\omega $ is a \textit{reflection}.
If $u_2\neq 0$, the isometry is called a \textit{glide reflection}. The
points in the ``mirror'' of reflection or glide reflection are $\alpha \in 
\mathbb{R}^2$ such that $\alpha _1={u_1}/2$.

If $\det \eta =-1$, one says that we have an \textit{opposite isometry}.
Every opposite isometry is a reflection or a glide reflection\cite{Armstrong}.

A \textit{planer image} is a function $\xi :\mathbb{R}^2\rightarrow D$, where
$D$ is a set, $D\neq \emptyset $.

From now on, let $\Omega _{\xi _{}}\equiv \Omega $ be the group of
isometries $\omega $ that leave $\xi $ invariant:
\begin{equation*}
\Omega _{\xi _{}}\equiv \Omega =\left\{ \omega \in \mathcal{I}:\xi \circ
\omega =\xi \right\} \text{.}
\end{equation*}

In the following, $\Omega ^{+}$ denotes the subgroup of $\Omega $ of the
isometries that preserve the orientation (direct isometries) and $\Omega
^{-} $ denotes the subgroup of $\Omega $ of the isometries that reverse the
orientation (opposite isometries).

We denote by $\Omega _{*}$ the group of orthogonal transformations
associated with the isometries of $\Omega $: 
\begin{equation*}
\Omega _{*}=\left\{ \eta \in O_2:\left( u,\eta \right) \in \Omega \right\}
\end{equation*}

In the following, if $\Lambda $ is a finite set, then $\left| \Lambda
\right| $ denotes its cardinal. Hence, if $G$ is a finite group, $\left|
G\right| $ denotes its order. For $a\in G$, $o\left( a\right) $ is the order
of the group generated by $a$ and $o\left( G\right) =\max \left\{ o\left(
a\right) :a\in G\right\} $.

\subsection{The seventeen wallpaper groups}

Let $\xi $ be a planer image, as before. From now on, we assume that $\Omega 
$ is a discrete group of isometries invariant under two linearly independent
translations which are of minimal length. Notice that $\Omega _{*}$ is
finite. We say that $\xi $ is a \textit{pattern} and that $\Omega $ is a 
\textit{wallpaper group}.

As it is very well known there are seventeen wallpaper groups. See \cite
{Armstrong}, \cite{Schattschneider}.

\begin{figure}[h]
  \centering
  \includegraphics[width=5.2702in]{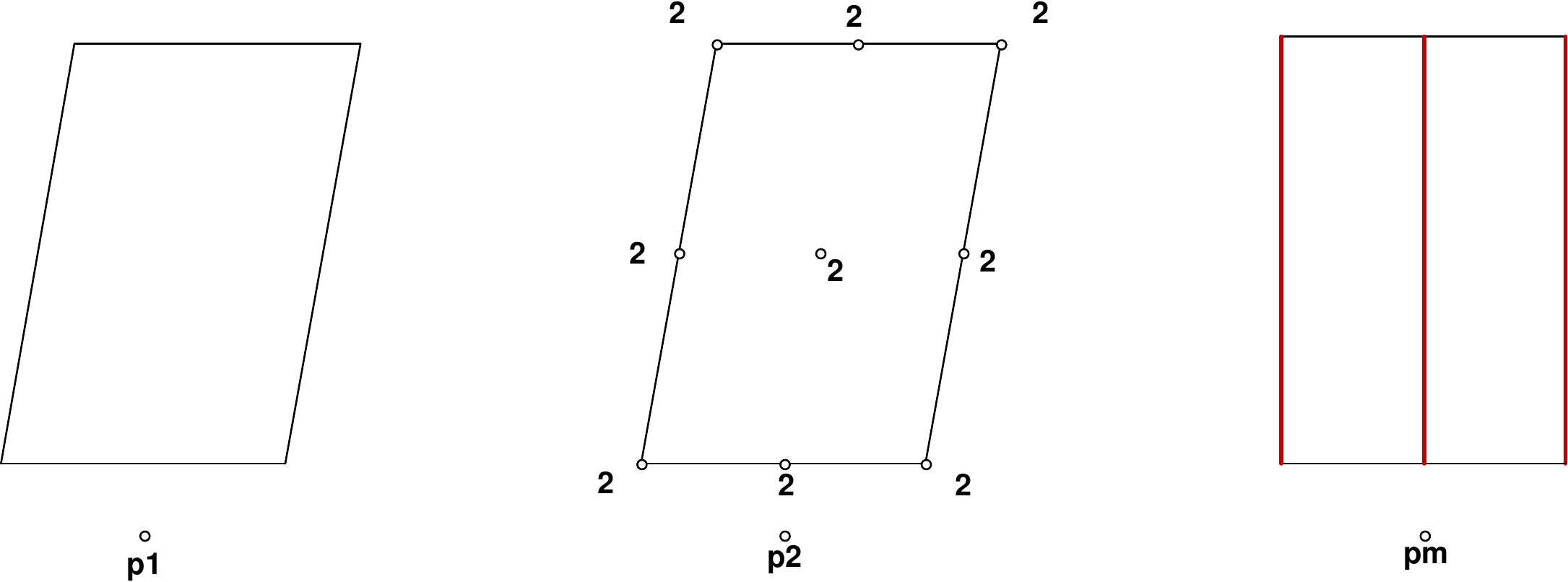}
  \caption{}
  \label{fig:02}
\end{figure}

\begin{figure}[h]
  \centering
  \includegraphics[width=5.271in]{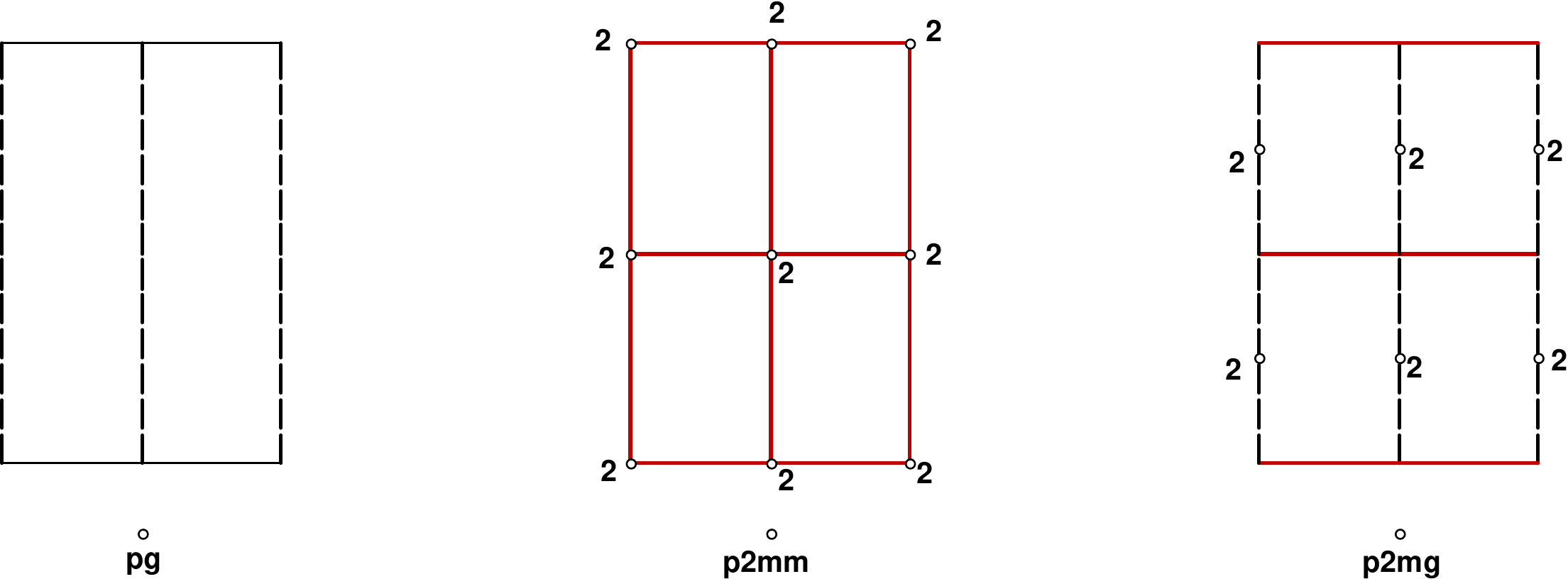}
  \caption{}
  \label{fig:03}
\end{figure}

\begin{figure}[H]
  \centering
  \includegraphics[width=5.2728in]{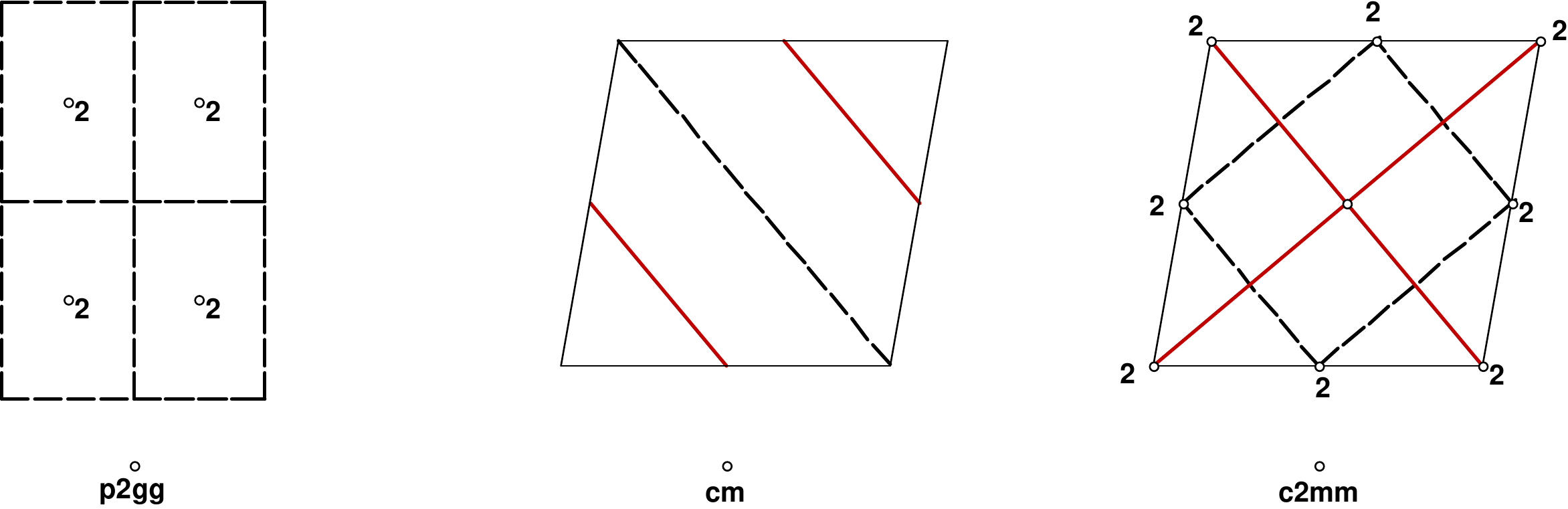}
  \caption{}
  \label{fig:04}
\end{figure}

\begin{figure}[H]
  \centering
  \includegraphics[width=5.2226in]{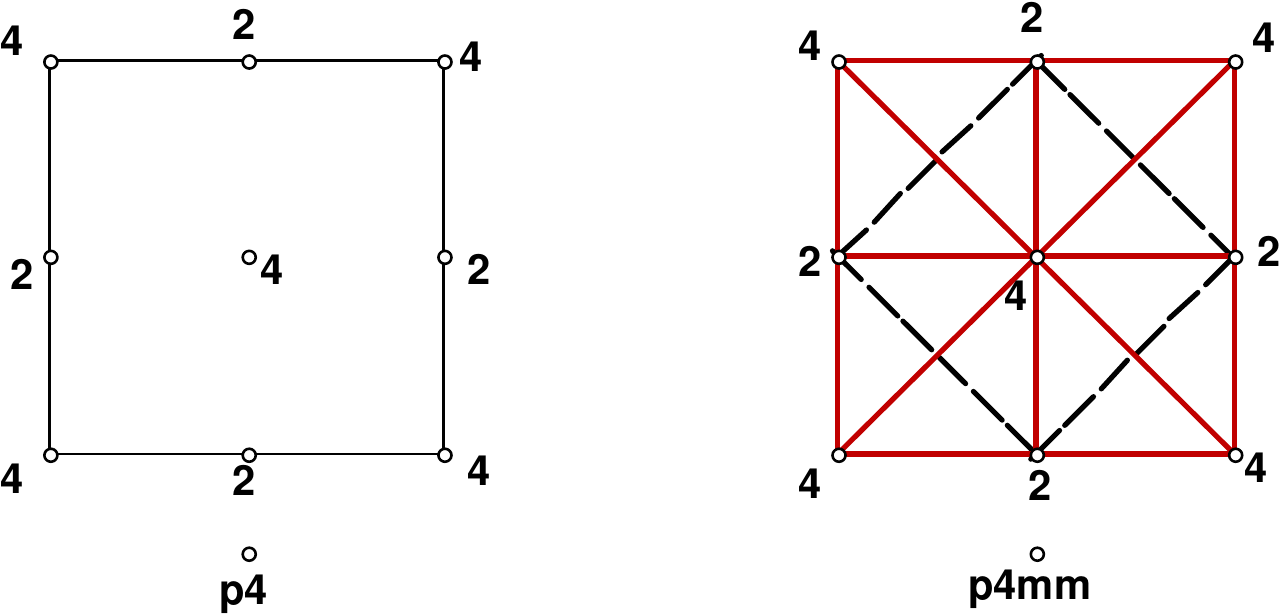}
  \caption{}
  \label{fig:05}
\end{figure}

\begin{figure}[H]
  \centering
  \includegraphics[width=4.95in]{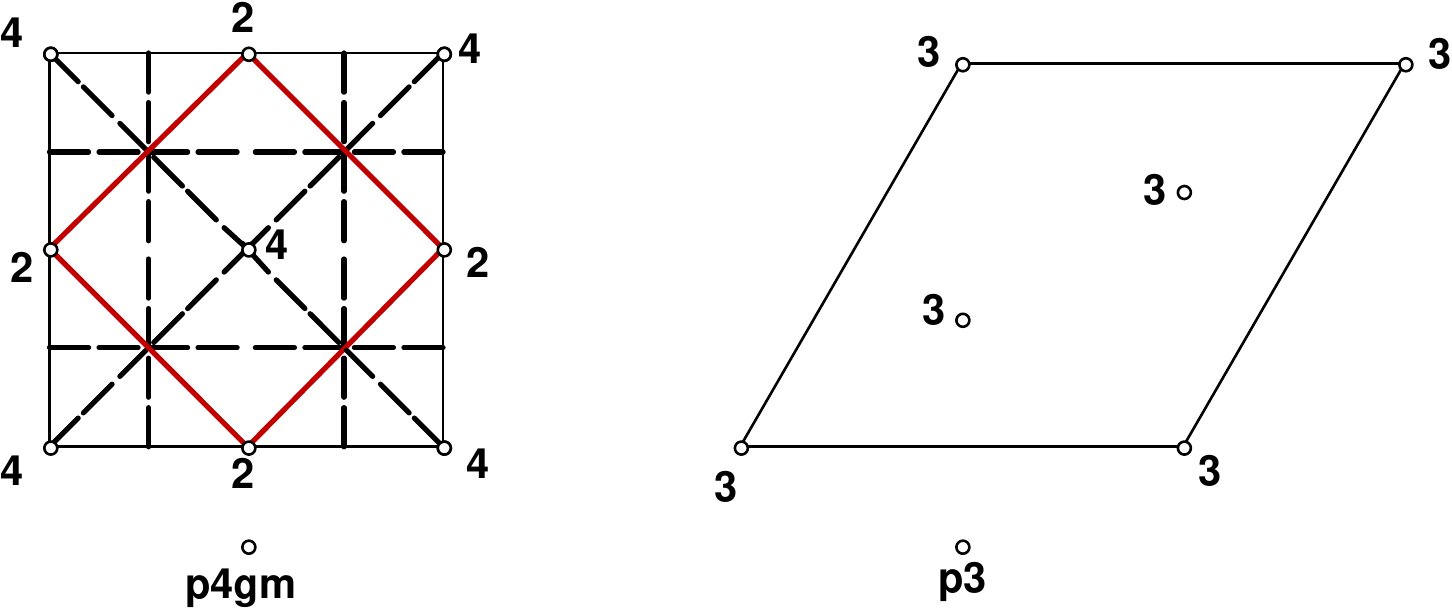}
  \caption{}
  \label{fig:06}
\end{figure}

\begin{figure}[H]
  \centering
  \includegraphics[width=5.78in]{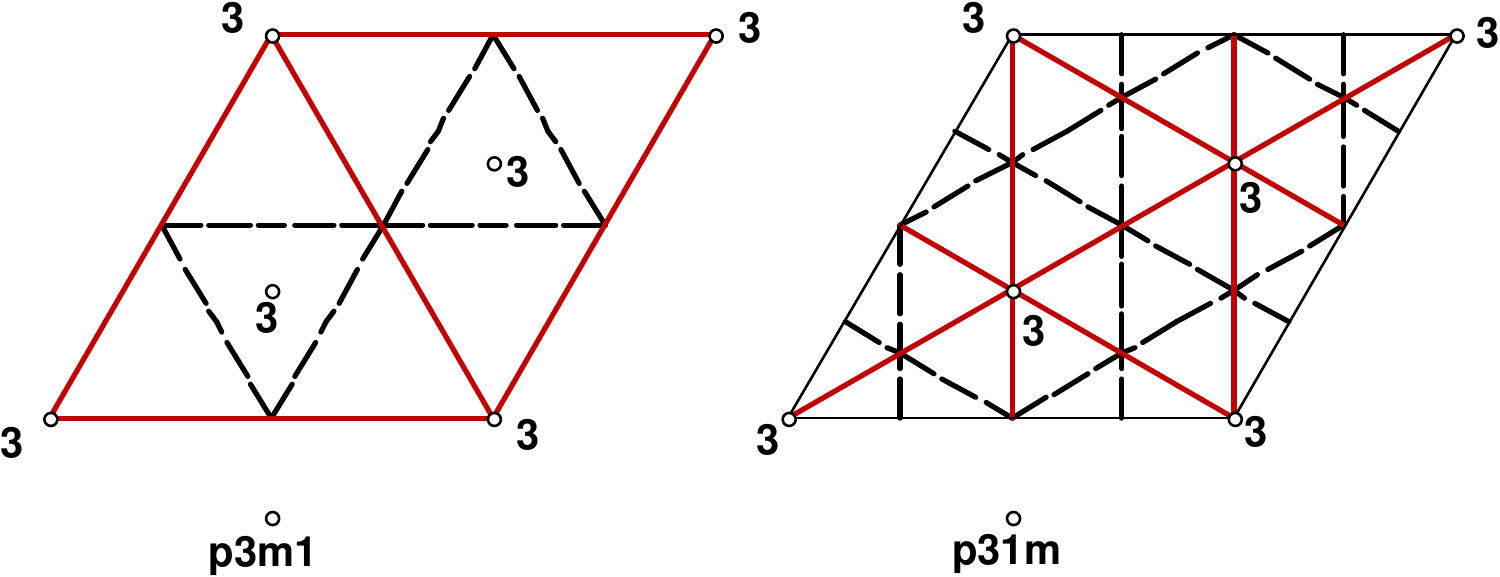}
  \caption{}
  \label{fig:07}
\end{figure}

\begin{figure}[H]
  \centering
  \includegraphics[width=5.78in]{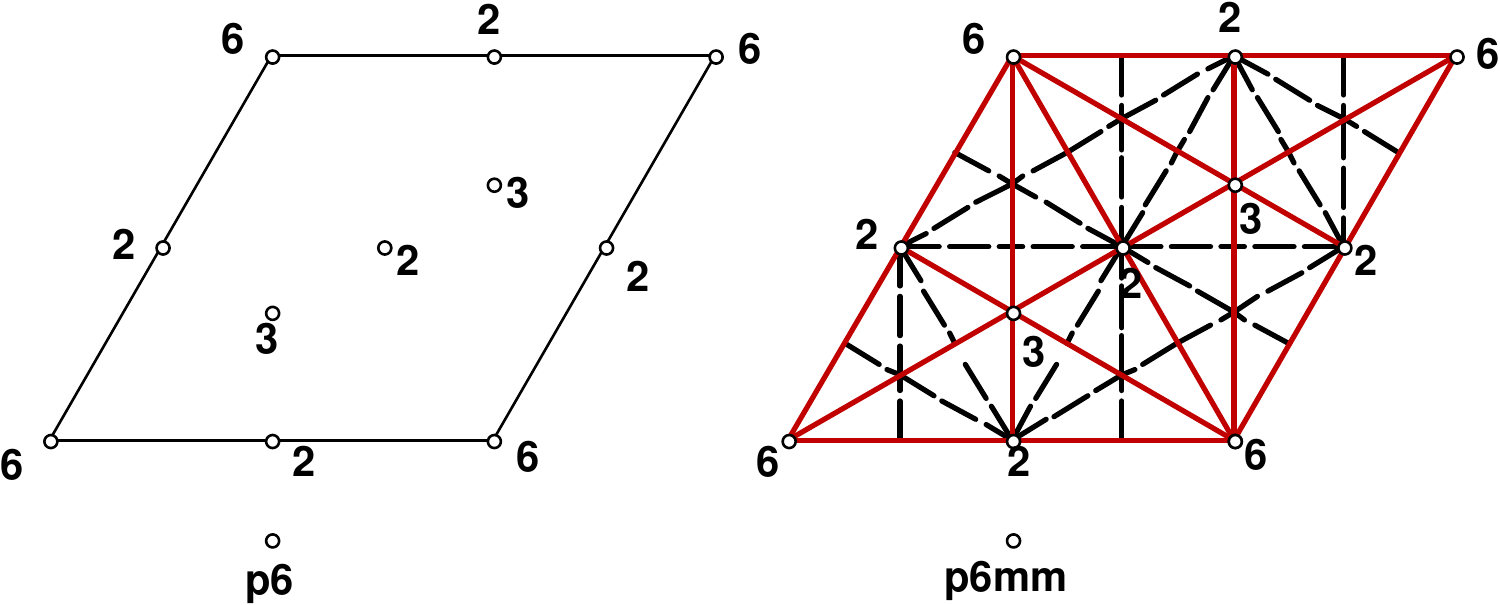}
  \caption{}
  \label{fig:08}
\end{figure}

Figures 2--8 represent fundamental regions of these seventeen groups. We use
the notations of Ref. \cite{Armstrong}. Here, small circles mean rotations
and numbers their order. Mirrors are represented by thick red lines and and
glides by broken ones.

\subsection{Wallpaper groups and permutation groups}

In this article, $S_n$ denotes the group of all permutations of $\left\{
1,2,\ldots ,n\right\} $; $a\in \ S_n$ means that $a$ is a one-to-one
function $a:\left\{ 1,2,\ldots ,n\right\} \rightarrow \left\{ 1,2,\ldots
,n\right\} $. The identity is $i$: $i\left( 1\right) =1,i\left( 2\right)
=2,\ldots ,i\left( n\right) =n$. If $a_1,a_2\in S_n$, we shall denote 
$a_1a_2\equiv a_1\circ a_2$.

We shall write $a=\left( m_1m_2\cdots m_k\right) \cdots \left( n_1n_2\cdots
n_l\right) $, if 
\begin{gather*}
a\left( m_1\right) =m_2,
\ a\left( m_2\right) =m_3,
\ \ldots,
\ a\left( m_k\right)=m_1,
\ \ldots , \\
a\left( n_1\right) =n_2,
\ a\left( n_2\right) =n_3,
\quad\ldots ,
\ a\left( n_k\right) =n_1
\end{gather*}
where $m_1,m_2,\ldots ,m_k,\ldots ,n_1,n_2,\ldots ,n_l\in \left\{ 1,2,\ldots
,n\right\} $.

If $p\in \left\{ 1,2,\ldots ,n\right\} \diagdown \left\{ m_1,m_2,\ldots
,m_k,\ldots ,n_1,n_2,\ldots ,n_l\right\} $, then $a\left( p\right) =p$.

The permutation $\left( m_1m_2\cdots m_k\right) $ is called a cyclic
permutation, or a cycle\ (in this case a $k$-cycle); $k$ is the length of
the cyclic permutation.

Let $S_n^{\pm }=\left\{ -1,1\right\} \times S_n$. If $\left( \delta
_1,a_1\right) ,\left( \delta _2,a_2\right) \in S_n^{\pm }$, then $\left(
\delta _1,a_1\right) \left( \delta _2,a_2\right) =\left( \delta _1\delta
_2,a_1a_2\right) $. $S_n^{\pm }$ is a group. We shall note $\left(
1,a\right) \equiv a$, $\left( -1,a\right) \equiv a_{-}$.

On permutation groups, see \cite{Wilson}.

From now on, let $\Omega $ be a wallpaper group and let $\zeta :\Omega
\rightarrow S_n^{\pm }$ be an homomorphism such that, if $\zeta \left(
\omega \right) =\left( \delta ,a\right) $, then $\delta =1$ if $\omega $
preserves the orientation and $\delta =-1$ otherwise (a reflection or a
glide reflection).

We shall say that the pair $\left( \Omega ,\zeta \right) $ is a
\textit{wallpaper group with permutations}. See \cite{Levine},
\cite{Schwarzenberger}.

In the following two sections we impose that $\left( \Omega ,\zeta \right) $
is connected, i.e., if $n_1,n_2\leq n$, then there exists $\omega \in \Omega 
$ such that $\zeta \left( \omega \right) \left( n_1\right) =n_2$.

\section{Translations}

Consider two independent vectors of the plane $\mathbb{R}^2$, $u$ and $v$, and
$\Omega $ the group of translations generated by them: $\Omega =\left\{
pu+qv:p,q\in \mathbb{Z}\right\} $. If $\zeta :\Omega \rightarrow S_n$ is a
group homomorphism, then $\zeta \left( u\right) $ and $\zeta \left( v\right) 
$ commute and, of course, $\zeta \left( pu+qv\right) =\zeta \left( u\right)
^p\zeta \left( v\right) ^q$. We assume that if $\zeta \left( \Omega \right)
\subset S_m$, then $m\geq n$.

One can identify $\Omega $ with $\mathbb{Z}^2$ by $\left( pu+qv\right)
\leftrightarrow \left( p,q\right) $ and define 
\begin{equation*}
\tilde{\zeta}\left( p,q\right) =\zeta \left( pu+qv\right) \left( 1\right) 
\text{.}
\end{equation*}

Notice that $\zeta $ is connected if and only if
$\tilde{\zeta}\left( \mathbb{Z}^2\right) =\left\{ 1,2,\ldots ,n\right\} $.

\begin{figure}[H]
  \centering
  \includegraphics[width=3.6703in]{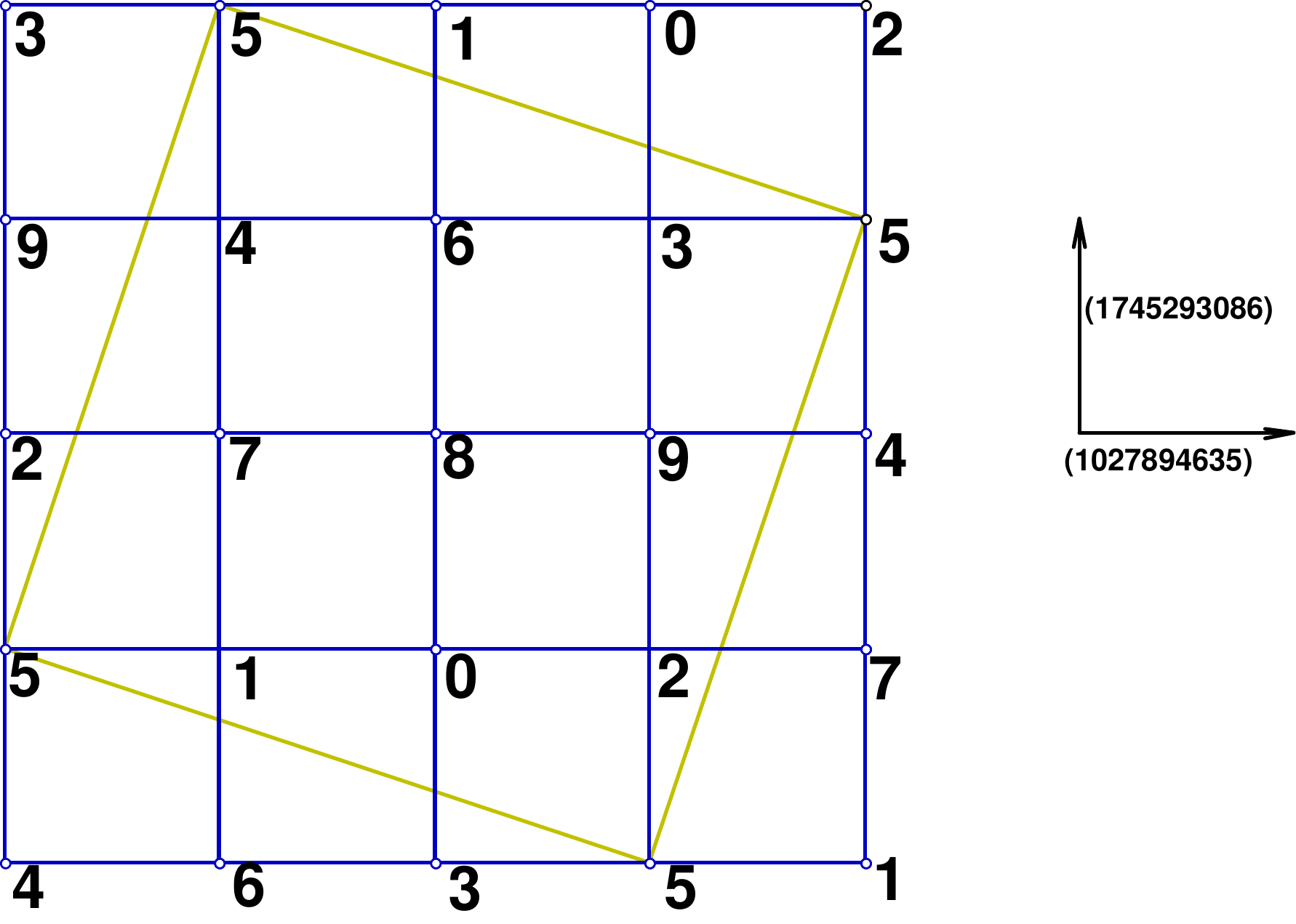}
  \caption{}
  \label{fig:09}
\end{figure}

As $\mathbb{Z}^2\subset \mathbb{R}^2$, $\tilde{\zeta}$ generates a periodical
pattern in the plane. A fundamental region of this pattern is a
parallelogram and it contains exactly $n$ points, (see Figures 9, 10 (a), 19
(a), 20 (a), 22 (a)). The vertices of this parallelogram in Figure 10 are
$\left( 0,0\right) $, $\left( p_1,q_1\right) $, $\left( p_2,q_2\right) $ and
$\left( p_1+p_2,q_1+q_2\right) $. Hence 
\begin{equation*}
n=\left| p_1q_2-p_2q_1\right|
\end{equation*}

The number of the $\zeta \left( u\right) $ cycles ($\mu _q$) and the order
of the $\zeta \left( u\right) $ cycles are 
\begin{equation*}
\mu _q=\gcd \left( \left| q_1\right| ,\left| q_2\right| \right)
\text{, }\frac n{\mu _q}\text{,}
\end{equation*}
and the number of the $\zeta \left( v\right) $ cycles ($\mu _p$) and the
order of the $\zeta \left( v\right) $ cycles are 
\begin{equation*}
\mu _p=\gcd \left( \left| p_1\right| ,\left| p_2\right| \right) 
\text{, }\frac n{\mu _p}\text{.}
\end{equation*}

\begin{figure}[H]
  \centering
  \includegraphics[width=5.7121in]{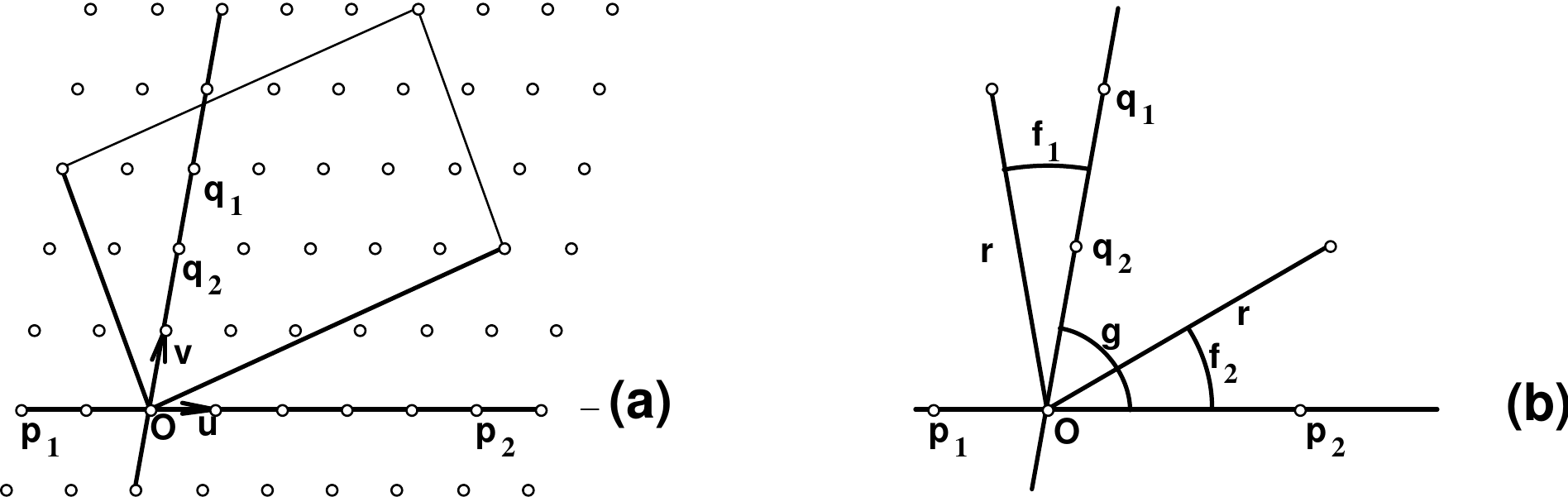}
  \caption{}
  \label{fig:10}
\end{figure}

Assume now that $\left\| u\right\| =\left\| v\right\| $ and that the
parallelogram is equilateral (see Figure 10 (b)). Then 
\begin{equation*}
\frac r{\sin g}=\frac{-p_1}{\sin f_1}=\frac{q_1}{\sin \left( g+f_1\right) }=
\frac{p_2}{\sin \left( g-f_2\right) }=\frac{q_2}{\sin f_2}\text{.}
\end{equation*}

1) If $f_1=f_2$, then 
\begin{equation*}
q_2=-p_1\text{,}\quad p_2=q_1+2p_1\cos g\text{,}
\end{equation*}
\begin{equation*}
\begin{split}
n &=p_1^2+q_1^2+2p_1q_1\cos g \\
&= p_2^2+q_2^2+2p_2q_2\cos g\text{.}
\end{split}
\end{equation*}

a) For $f_1=f_2$, and $g=\frac \pi 3$ (a regular triangular grid), then 
\begin{equation*}
q_2=-p_1\text{,}\quad
p_2=p_1+q_1\text{,}\quad
n=p_1^2+q_1^2+p_1q_1=p_2^2+q_2^2+p_2q_2\text{.}
\end{equation*}

The next table presents values of $n=p^2+q^2+pq$.
\begin{equation*}
\begin{tabular}{|c|c|c|c|c|c|c|c|c|}
\hline
$
\begin{array}{l}
p\rightarrow \\ 
q\downarrow
\end{array}
$ & $0$ & $1$ & $2$ & $3$ & $4$ & $5$ & $6$ & $\cdots $ \\ \hline
$0$ & $0$ & $1$ & $4$ & $9$ & $16$ & $25$ & $36$ & $\cdots $ \\ \hline
$1$ & $1$ & $3$ & $7$ & $13$ & $21$ & $31$ & $43$ & $\cdots $ \\ \hline
$2$ & $4$ & $7$ & $12$ & $19$ & $28$ & $39$ & $52$ & $\cdots $ \\ \hline
$3$ & $9$ & $13$ & $19$ & $27$ & $37$ & $49$ & $63$ & $\cdots $ \\ \hline
$4$ & $16$ & $21$ & $28$ & $37$ & $48$ & $61$ & $76$ & $\cdots $ \\ \hline
$\cdots $ & $\cdots $ & $\cdots $ & $\cdots $ & $\cdots $ & $\cdots $ & $\cdots $ & $\cdots $ & $\cdots $ \\ \hline
\end{tabular}
\end{equation*}
Notice that the number of equilateral triangles in the fundamental region is 
$2n$.

b) For $f_1=f_2$, and $g=\frac \pi 4$ (a square grid), then 
\begin{equation*}
q_2=-p_1\text{,}\quad
p_2=q_1\text{,}\quad
n=p_1^2+q_1^2=p_2^2+q_2^2\text{.}
\end{equation*}

The next table presents values of $n=p^2+q^2$.
\begin{equation*}
\begin{tabular}{|c|c|c|c|c|c|c|c|c|}
\hline
$
\begin{array}{l}
p\rightarrow \\ 
q\downarrow
\end{array}
$ & $0$ & $1$ & $2$ & $3$ & $4$ & $5$ & $6$ & $\cdots $ \\ \hline
$0$ & $0$ & $1$ & $4$ & $9$ & $16$ & $25$ & $36$ & $\cdots $ \\ \hline
$1$ & $1$ & $2$ & $5$ & $10$ & $17$ & $26$ & $37$ & $\cdots $ \\ \hline
$2$ & $4$ & $5$ & $8$ & $13$ & $20$ & $29$ & $40$ & $\cdots $ \\ \hline
$3$ & $9$ & $10$ & $13$ & $18$ & $25$ & $34$ & $45$ & $\cdots $ \\ \hline
$4$ & $16$ & $17$ & $20$ & $25$ & $32$ & $41$ & $52$ & $\cdots $ \\ \hline
$\cdots $ & $\cdots $ & $\cdots $ & $\cdots $ & $\cdots $ & $\cdots $ & $\cdots $ & $\cdots $ & $\cdots $ \\ \hline
\end{tabular}
\end{equation*}

2) If $f_1=-f_2$, then 
\begin{equation*}
q_2=p_1\text{,}\quad
p_2=q_1\text{,}\quad
n=\left| p_1^2-q_1^2\right| =\left|
p_2^2-q_2^2\right| \text{.}
\end{equation*}

Examples: 
\begin{equation*}
\begin{tabular}{|c|c|c|c|c|c|c|c|c|c|}
\hline
& $p_1$ & $q_1$ & $p_2$ & $q_2$ & $n$ & $\mu _q$ & $\frac n{\mu _q}$ & $\mu
_p$ & $\frac n{\mu _p}$ \\ \hline
Figure 9 & $-3$ & $1$ & $1$ & $3$ & $10$ & $1$ & $10$ & $1$ & $10$ \\ \hline
Figure 10 (a) & $-2$ & $3$ & $5$ & $2$ & $19$ & $1$ & $19$ & $1$ & $19$ \\ 
\hline
Figure 19 (a) & $-1$ & $2$ & $2$ & $1$ & $5$ & $1$ & $5$ & $1$ & $5$ \\ 
\hline
Figure 20 (a) & $-2$ & $2$ & $2$ & $2$ & $8$ & $2$ & $4$ & $2$ & $4$ \\ 
\hline
Figure 22 (a) & $-3$ & $1$ & $-2$ & $3$ & $7$ & $1$ & $7$ & $1$ & $7$ \\ 
\hline
\end{tabular}
\end{equation*}

\section{Rotations and reflections}

In this section we describe the different possibilities for wallpaper groups
with permutations. Here, the translations are generated by rotations,
reflections and glide reflections. We impose that $\left( \Omega ,\zeta
\right) $ is connected and that $n\leq o\left( \Omega _{*}\right) $. As
before, in the figures, small circles mean rotations and numbers their
order. Mirrors are represented by thick red lines and and glides by broken
ones. From now on, the letters $a$, $b$, $c$, $d$, $x$, $y$, $z$, $w$, $u$
and $v$ represent permutations.

Let us describe, briefly, the method we follow in this section.

If $\omega _1$ is a rotation of order $k$, and $\omega _2$ is a rotation of
order $j$, then $\omega _1$ transforms $\omega _2$ in another rotation of
order $j$, $\omega _3$, which is $\omega _1\omega _2\omega _1^{-1}\equiv
\omega _1\left( \omega _2\right) $.

If one wants to translate the isometries into elements of $S_n^{\pm }$ the
function must be such that if $\omega _1\longmapsto s_1$, $\omega
_2\longmapsto s_2$, then $\omega _3\equiv \omega _1\left( \omega _2\right)
\longmapsto s_1s_2s_1^{-1}$.

Note that if $s_1=\left( \delta _1,a_1\right) $, $s_2=\left( \delta
_2,a_2\right) $, where 
\begin{equation*}
a_2=\left( m_1m_2\cdots m_l\right) \cdots \text{,}
\end{equation*}
then $\omega _3\equiv \omega _1\left( \omega _2\right) \longmapsto
s_3=\left( \delta _2,a_3\right) $, with 
\begin{equation*}
a_3=\left( a_1\left( m_1\right) a_1\left( m_2\right) \cdots a_1\left(
m_l\right) \right) \cdots \text{.}
\end{equation*}

We assign to every axis of order $k$\ ($\equiv \omega $), a permutation $a$
of order $k_1$, a divisor of $k$, so that to the counter clock-wise rotation
of $\frac{2\pi }k$, $\omega $, corresponds the permutation $a$. This
association must be coherent in the sense that it generates a group
homomorphism.

If we assign to $\omega $\ and $\omega _1$\ (two neighbor axes) the
permutations $a$ and $a_1$, then to the axis $\omega \left( \omega _1\right)
=\omega \omega _1\omega ^{-1}$ we must assign $aa_1a^{-1}$. When $a_1=\left(
m_1m_2\cdots m_l\right) \cdots $, then $aa_1a^{-1}=\left( a\left( m_1\right)
a\left( m_2\right) \cdots a\left( m_l\right) \right) \cdots $.

\begin{figure}[H]
  \centering
  \includegraphics[width=5.271in]{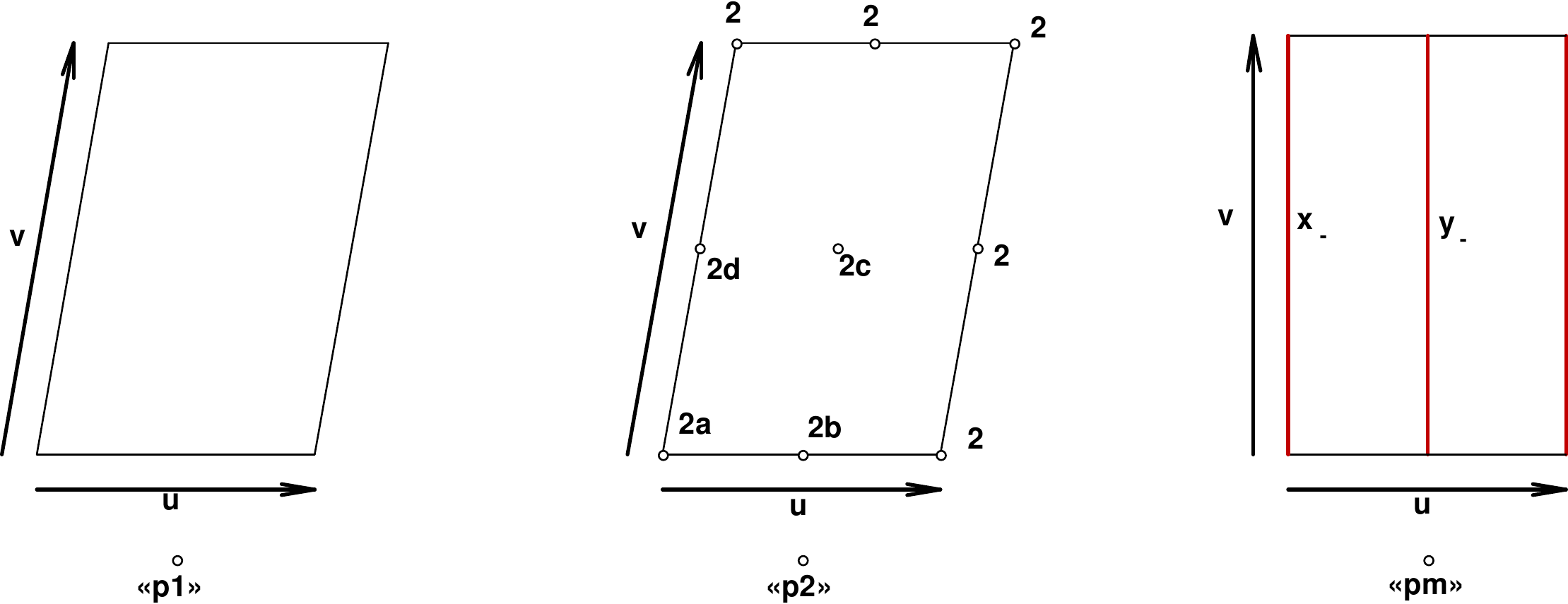}
  \caption{}
  \label{fig:11}
\end{figure}

\subsection{\boldmath Figure 11, \guillemotleft $p2$\guillemotright}
In the case \guillemotleft $p2$\guillemotright , the possibilities for $a$,
$b$, $c$ and $d$ are $i$ or $\left( 12\right) $. Notice that $c=dba$, $u=ba$
and $v=da$.
\begin{equation*}
\begin{tabular}{|c|c|c|c|c|c|}
\hline
$a$ & $b$ & $d$ & $c$ & $u$ & $v$ \\ \hline
$i$ & $i$ & $i$ & $i$ & $i$ & $i$ \\ \hline
$\left( 12\right) $ & $i$ & $i$ & $\left( 12\right) $ & $\left( 12\right) $
& $\left( 12\right) $ \\ \hline
$\left( 12\right) $ & $i$ & $\left( 12\right) $ & $i$ & $\left( 12\right) $
& $i$ \\ \hline
$\left( 12\right) $ & $\left( 12\right) $ & $\left( 12\right) $ & $\left(
12\right) $ & $i$ & $i$ \\ \hline
\end{tabular}
\end{equation*}

\subsection{\boldmath Figure 11, \guillemotleft $pm$\guillemotright}

In the case \guillemotleft $pm$\guillemotright , the possibilities for $x$
and $y$ are $i$ or $\left( 12\right) $. Notice that $u=xy$ and $v=i$.
\begin{equation*}
\begin{tabular}{|c|c|c|c|}
\hline
$x$ & $y$ & $u$ & $v$ \\ \hline
$i$ & $i$ & $i$ & $i$ \\ \hline
$\left( 12\right) $ & $i$ & $\left( 12\right) $ & $i$ \\ \hline
$\left( 12\right) $ & $\left( 12\right) $ & $i$ & $i$ \\ \hline
\end{tabular}
\end{equation*}

\subsection{\boldmath Figure 12, \guillemotleft $pg$\guillemotright}

In the case \guillemotleft $pg$\guillemotright , the possibilities for $x$
and $y$ are $i$ or $\left( 12\right) $. Notice that $u=xy$ and $v=i$.
\begin{equation*}
\begin{tabular}{|c|c|c|c|}
\hline
$x$ & $y$ & $u$ & $v$ \\ \hline
$i$ & $i$ & $i$ & $i$ \\ \hline
$\left( 12\right) $ & $i$ & $\left( 12\right) $ & $i$ \\ \hline
$\left( 12\right) $ & $\left( 12\right) $ & $i$ & $i$ \\ \hline
\end{tabular}
\end{equation*}
\begin{figure}[H]
  \centering
  \includegraphics[width=5.271in]{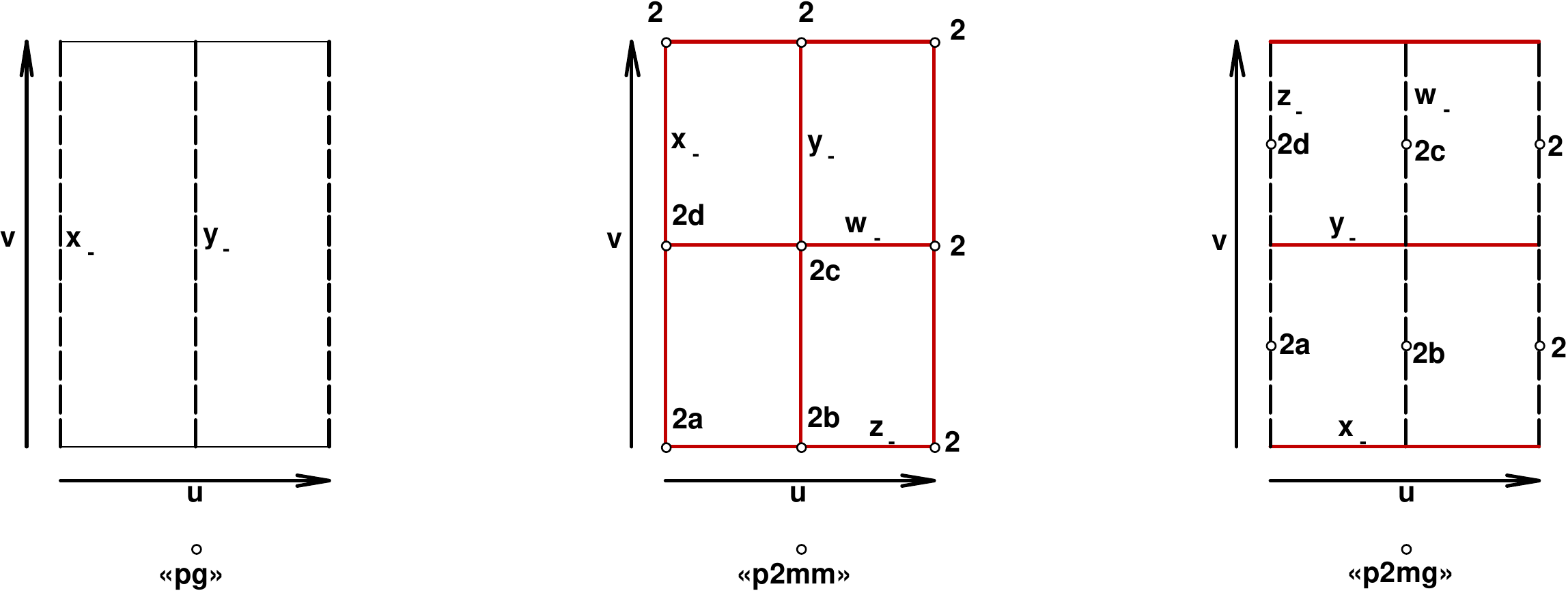}
  \caption{}
  \label{fig:12}
\end{figure}

\subsection{\boldmath Figure 12, \guillemotleft $p2mm$\guillemotright}

In the case \guillemotleft $p2mm$\guillemotright , the possibilities for
$a$, $b$, $c$ and $d$ are $i$ or $\left( 12\right) $ as in \guillemotleft $p2$\guillemotright . The possibilities for $x$ and $y$ are $i$ or $\left(
12\right) $.

Notice that $c=dba$, $u=ba$, $v=da$, $y=ux$, $z=ax$, and $w=dx$.
\begin{equation*}
\begin{tabular}{|c|c|c|c|c|c|c|c|c|c|}
\hline
$a$ & $b$ & $d$ & $c$ & $u$ & $v$ & $x$ & $y$ & $z$ & $w$ \\ \hline
$i$ & $i$ & $i$ & $i$ & $i$ & $i$ & $i$ & $i$ & $i$ & $i$ \\ \hline
$i$ & $i$ & $i$ & $i$ & $i$ & $i$ & $\left( 12\right) $ & $\left( 12\right) $ & $\left( 12\right) $ & $\left( 12\right) $ \\ \hline
$\left( 12\right) $ & $i$ & $i$ & $\left( 12\right) $ & $\left( 12\right) $ & $\left( 12\right) $ & $i$ & $\left( 12\right) $ & $\left( 12\right) $ & $i$ \\ \hline
$\left( 12\right) $ & $i$ & $\left( 12\right) $ & $i$ & $\left( 12\right) $
& $i$ & $i$ & $\left( 12\right) $ & $\left( 12\right) $ & $\left( 12\right) $
\\ \hline
$\left( 12\right) $ & $i$ & $\left( 12\right) $ & $i$ & $\left( 12\right) $
& $i$ & $\left( 12\right) $ & $i$ & $i$ & $i$ \\ \hline
$\left( 12\right) $ & $\left( 12\right) $ & $\left( 12\right) $ & $\left(
12\right) $ & $i$ & $i$ & $i$ & $i$ & $\left( 12\right) $ & $\left(
12\right) $ \\ \hline
\end{tabular}
\end{equation*}

\subsection{\boldmath Figure 12, \guillemotleft $p2mg$\guillemotright}

In the case \guillemotleft $p2mg$\guillemotright , the possibilities for $a$,
$b$, $c$ and $d$ are $i$ or $\left( 12\right) $ as in
\guillemotleft $p2$\guillemotright . The possibilities for $x$ and $y$
are $i$ or $\left(12\right) $.

Notice that $d=a$, $c=b$, $u=ba$, $v=i$, $y=x$, $z=ax$, $w=bx$.
\begin{equation*}
\begin{tabular}{|c|c|c|c|c|c|c|c|c|}
\hline
$a$ & $b$ & $d$ & $c$ & $u$ & $v$ & $x,y$ & $z$ & $w$ \\ \hline
$i$ & $i$ & $i$ & $i$ & $i$ & $i$ & $i$ & $i$ & $i$ \\ \hline
$i$ & $i$ & $i$ & $i$ & $i$ & $i$ & $\left( 12\right) $ & $\left( 12\right) $
& $\left( 12\right) $ \\ \hline
$\left( 12\right) $ & $i$ & $\left( 12\right) $ & $i$ & $\left( 12\right) $
& $i$ & $i$ & $\left( 12\right) $ & $i$ \\ \hline
$\left( 12\right) $ & $i$ & $\left( 12\right) $ & $i$ & $\left( 12\right) $
& $i$ & $\left( 12\right) $ & $i$ & $\left( 12\right) $ \\ \hline
$\left( 12\right) $ & $\left( 12\right) $ & $\left( 12\right) $ & $\left(
12\right) $ & $i$ & $i$ & $i$ & $\left( 12\right) $ & $\left( 12\right) $ \\ 
\hline
$\left( 12\right) $ & $\left( 12\right) $ & $\left( 12\right) $ & $\left(
12\right) $ & $i$ & $i$ & $\left( 12\right) $ & $i$ & $i$ \\ \hline
\end{tabular}
\end{equation*}

\subsection{\boldmath Figure 13, \guillemotleft $p2gg$\guillemotright}

In the case \guillemotleft $p2gg$\guillemotright , the possibilities for
$a$, $b$, $c$ and $d$ are $i$ or $\left( 12\right) $ as in
\guillemotleft $p2$\guillemotright . The possibilities for $x$ and $y$ are
$i$ or $\left(12\right) $.

Notice that $a=b=c=d$, $u=i$, $v=i$, $y=x$, $z=w=ax$.
\begin{equation*}
\begin{tabular}{|c|c|c|c|c|c|c|c|}
\hline
$a$ & $b$ & $d$ & $c$ & $u$ & $v$ & $x,y$ & $z,w$ \\ \hline
$i$ & $i$ & $i$ & $i$ & $i$ & $i$ & $i$ & $i$ \\ \hline
$i$ & $i$ & $i$ & $i$ & $i$ & $i$ & $\left( 12\right) $ & $\left( 12\right) $
\\ \hline
$\left( 12\right) $ & $\left( 12\right) $ & $\left( 12\right) $ & $\left(
12\right) $ & $i$ & $i$ & $i$ & $\left( 12\right) $ \\ \hline
\end{tabular}
\end{equation*}

\begin{figure}[H]
  \centering
  \includegraphics[width=5.2728in]{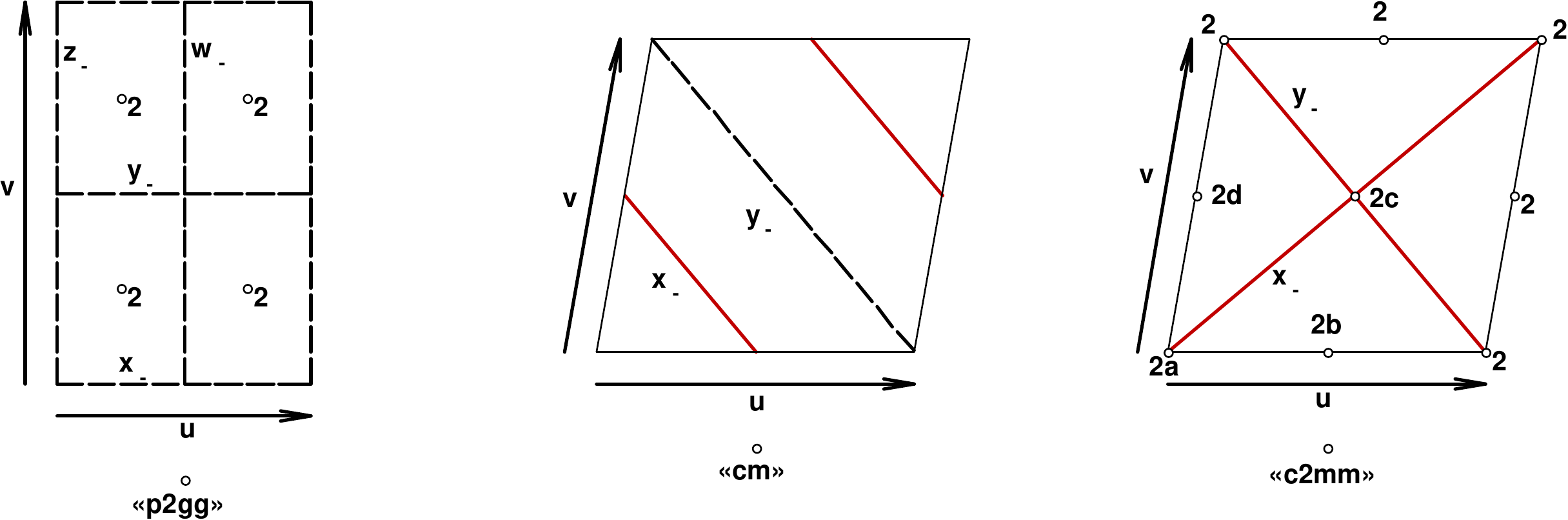}
  \caption{}
  \label{fig:13}
\end{figure}

\subsection{\boldmath Figure 13, \guillemotleft $cm$\guillemotright}

In the case \guillemotleft $cm$\guillemotright , the possibilities for $x$
and $y$ are $i$ or $\left( 12\right) $. Notice that $u=v=xy$.
\begin{equation*}
\begin{tabular}{|c|c|c|c|}
\hline
$x$ & $y$ & $u$ & $v$ \\ \hline
$i$ & $i$ & $i$ & $i$ \\ \hline
$i$ & $\left( 12\right) $ & $\left( 12\right) $ & $\left( 12\right) $ \\ 
\hline
$\left( 12\right) $ & $i$ & $\left( 12\right) $ & $\left( 12\right) $ \\ 
\hline
$\left( 12\right) $ & $\left( 12\right) $ & $i$ & $i$ \\ \hline
\end{tabular}
\end{equation*}

\subsection{\boldmath Figure 13, \guillemotleft $c2mm$\guillemotright}

In the case \guillemotleft $c2mm$\guillemotright , the possibilities for
$a$, $b$, $c$ and $d$ are $i$ or $\left( 12\right) $ as in
\guillemotleft $p2$\guillemotright . The possibilities for $x$ and $y$ are
$i$ or $\left(12\right) $.

Notice that $b=d$, $c=a$, $u=v=ba$, $v=da$, $y=ax$.
\begin{equation*}
\begin{tabular}{|c|c|c|c|c|c|c|}
\hline
$a$ & $b$ & $d$ & $c$ & $u,v$ & $x$ & $y$ \\ \hline
$i$ & $i$ & $i$ & $i$ & $i$ & $i$ & $i$ \\ \hline
$i$ & $i$ & $i$ & $i$ & $i$ & $\left( 12\right) $ & $\left( 12\right) $ \\ 
\hline
$i$ & $\left( 12\right) $ & $\left( 12\right) $ & $i$ & $\left( 12\right) $
& $i$ & $i$ \\ \hline
$i$ & $\left( 12\right) $ & $\left( 12\right) $ & $i$ & $\left( 12\right) $
& $\left( 12\right) $ & $i$ \\ \hline
$\left( 12\right) $ & $i$ & $i$ & $\left( 12\right) $ & $\left( 12\right) $
& $i$ & $\left( 12\right) $ \\ \hline
$\left( 12\right) $ & $\left( 12\right) $ & $\left( 12\right) $ & $\left(12\right) $ & $i$ & $i$ & $\left( 12\right) $ \\ \hline
\end{tabular}
\end{equation*}

\subsection{\boldmath Figure 14, \guillemotleft $p4$\guillemotright}

In the case \guillemotleft $p4$\guillemotright , $b=ca$, $d=ac$, $u=ca^{-1}$,
$v=c^{-1}a$. We shall take always $a$ and $c$ such that $o\left( c\right)
\leq o\left( a\right) $.

The possibilities for $a$ are $i$ (for $n=1$) and permutations of the type
$\left( 12\right) $ (for $n=2$), $\left( 12\right) \left( 34\right) $ and
$\left( 1234\right) $ (for $n=4$).

If $a=i$, then $c=i$. If $a=\left( 12\right) $, then $c=i$ or $c=\left(
12\right) $.

If $a=\left( 12\right) \left( 34\right) $, then $c=\left( 13\right) \left(
24\right) $.

If $a=\left( 1234\right) $, then $cac^{-1}=\left( 1234\right) $ or $\left(
1432\right) $, hence the possibilities for $c$ are $\left( 1234\right) $,
$\left( 1432\right) $, $\left( 13\right) $, $\left( 24\right) $, $\left(
12\right) \left( 34\right) $ and $\left( 14\right) \left( 23\right) $.
Notice that, since $a\left( 13\right) a^{-1}=\left( 24\right) $ and $a\left(
12\right) \left( 34\right) a^{-1}=\left( 14\right) \left( 23\right) $, the
possibilities for $c$ are, in fact, $\left( 1234\right) $, $\left(
1432\right) $, $\left( 13\right) $ and $\left( 12\right) \left( 34\right) $.
\begin{equation*}
\begin{tabular}{|c|c|c|c|c|}
\hline
$a$ & $c$ & $b$ & $u$ & $v$ \\ \hline
$i$ & $i$ & $i$ & $i$ & $i$ \\ \hline
$\left( 12\right) $ & $i$ & $\left( 12\right) $ & $\left( 12\right) $ &
$\left( 12\right) $ \\ \hline
$\left( 12\right) $ & $\left( 12\right) $ & $i$ & $i$ & $i$ \\ \hline
$\left( 12\right) \left( 34\right) $ & $\left( 13\right) \left( 24\right) $
& $\left( 14\right) \left( 23\right) $ & $\left( 14\right) \left( 23\right) $
& $\left( 14\right) \left( 23\right) $ \\ \hline
$\left( 1234\right) $ & $\left( 1234\right) $ & $\left( 13\right) \left(
24\right) $ & $i$ & $i$ \\ \hline
$\left( 1234\right) $ & $\left( 1432\right) $ & $i$ & $\left( 13\right)
\left( 24\right) $ & $\left( 13\right) \left( 24\right) $ \\ \hline
$\left( 1234\right) $ & $\left( 13\right) $ & $\left( 12\right) \left(
34\right) $ & $\left( 14\right) \left( 23\right) $ & $\left( 12\right)
\left( 34\right) $ \\ \hline
$\left( 1234\right) $ & $\left( 12\right) \left( 34\right) $ & $\left(
24\right) $ & $\left( 13\right) $ & $\left( 24\right) $ \\ \hline
\end{tabular}
\end{equation*}

\begin{figure}[ht]
  \centering
  \includegraphics[width=4.03in]{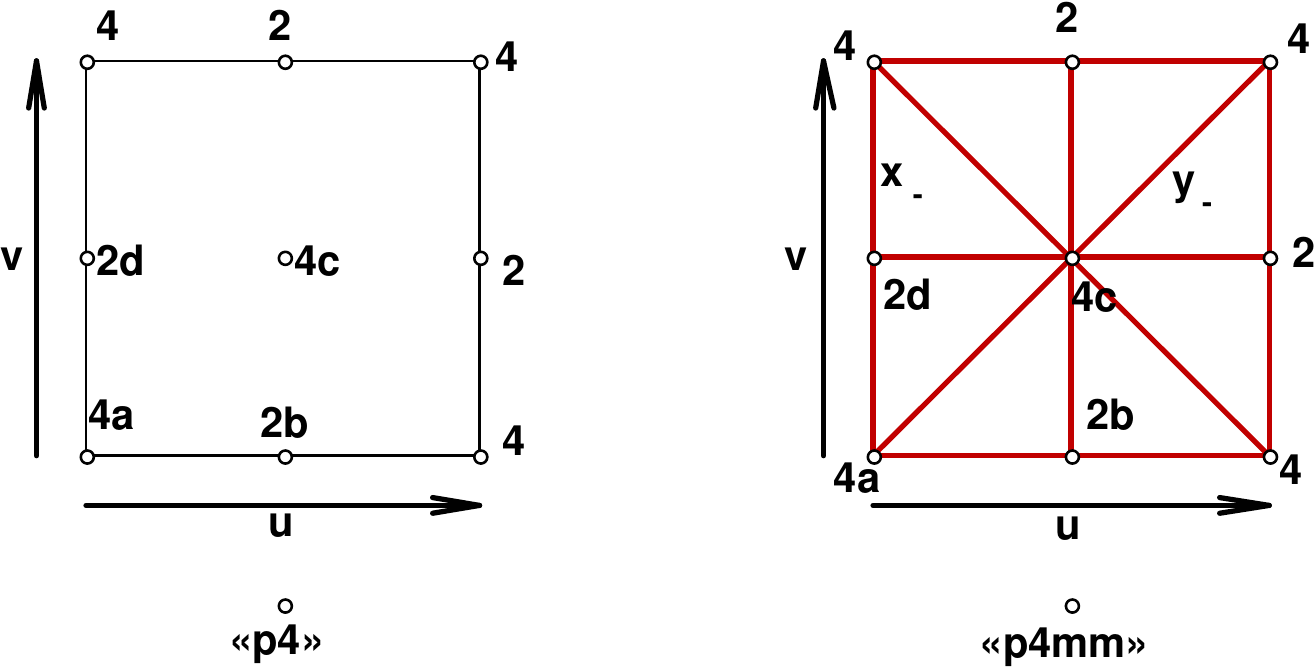}
  \vspace{-20pt}
  \caption{}
  \label{fig:14}
\end{figure}

\subsection{\boldmath Figure 14, \guillemotleft $p4mm$\guillemotright}

In the case \guillemotleft $p4mm$\guillemotright , $a$, $b$, $c$, $d$, $u$
and $v$ are as in \guillemotleft $p4$\guillemotright .

Here, $y=a^{-1}x$.

For $a=i$, $n\leq 2$, the possibilities for $x$ are $i$ and $\left(
12\right) $.

For $a=\left( 12\right) $, $n=2$, the possibilities for $x$ are $i$ and
$\left( 12\right) $.

For $a=\left( 12\right) \left( 34\right) $, $c=\left( 13\right) \left(
24\right) $, the possibilities for $x$ are $i$, $\left( 12\right) \left(
34\right) $, $\left( 13\right) \left( 24\right) $ and $\left( 14\right)
\left( 23\right) $.

For $a=\left( 1234\right) $, $c=\left( 1234\right) $, the possibilities for
$x$ are $\left( 13\right) $ and $\left( 12\right) \left( 34\right) $.

For $a=\left( 1234\right) $, $c=\left( 1432\right) $, the possibilities for
$x$ are $\left( 13\right) $ and $\left( 12\right) \left( 34\right) $.

For $a=\left( 1234\right) $, $c=\left( 13\right) $, the possibilities for $x$
are $\left( 12\right) \left( 34\right) $ and $\left( 14\right) \left(
23\right) $.

For $a=\left( 1234\right) $, $c=\left( 12\right) \left( 34\right) $, the
possibilities for $x$ are $\left( 13\right) $ and $\left( 24\right) $.

As \guillemotleft $p4mm$\guillemotright\ is constructed from
\guillemotleft $p4$\guillemotright\ adding reflections, these reflections can connect
permutations:

a) $a=\left( 12\right) \left( 34\right) $, $c=i$, $b=\left( 12\right) \left(
34\right) $, $u=\left( 12\right) \left( 34\right) $, $x=\left( 13\right)
\left( 24\right) $.

b) $a=\left( 12\right) \left( 34\right) $, $c=\left( 12\right) \left(
34\right) $ , $b=i$, $u=i$, $x=\left( 13\right) \left( 24\right) $.

\begin{equation*}
\begin{tabular}{|c|c|c|c|c|c|}
\hline
$a$ & $c$ & $b$ & $u$ & $x$ & $y$ \\ \hline
$i$ & $i$ & $i$ & $i$ & $i$ & $i$ \\ \hline
$i$ & $i$ & $i$ & $i$ & $\left( 12\right) $ & $\left( 12\right) $ \\ \hline
$\left( 12\right) $ & $i$ & $\left( 12\right) $ & $\left( 12\right) $ & $i$
& $\left( 12\right) $ \\ \hline
$\left( 12\right) $ & $i$ & $\left( 12\right) $ & $\left( 12\right) $ & $\left( 12\right) $ & $i$ \\ \hline
$\left( 12\right) \left( 34\right) $ & $i$ & $\left( 12\right) \left(
34\right) $ & $\left( 12\right) \left( 34\right) $ & $\left( 13\right)
\left( 24\right) $ & $\left( 14\right) \left( 23\right) $ \\ \hline
$\left( 12\right) $ & $\left( 12\right) $ & $i$ & $i$ & $i$ & $\left(
12\right) $ \\ \hline
$\left( 12\right) $ & $\left( 12\right) $ & $i$ & $i$ & $\left( 12\right) $
& $i$ \\ \hline
$\left( 12\right) \left( 34\right) $ & $\left( 12\right) \left( 34\right) $
& $i$ & $i$ & $\left( 13\right) \left( 24\right) $ & $\left( 14\right)
\left( 23\right) $ \\ \hline
$\left( 12\right) \left( 34\right) $ & $\left( 13\right) \left( 24\right) $
& $\left( 14\right) \left( 23\right) $ & $\left( 14\right) \left( 23\right) $
& $i$ & $\left( 12\right) \left( 34\right) $ \\ \hline
$\left( 12\right) \left( 34\right) $ & $\left( 13\right) \left( 24\right) $
& $\left( 14\right) \left( 23\right) $ & $\left( 14\right) \left( 23\right) $
& $\left( 12\right) \left( 34\right) $ & $i$ \\ \hline
$\left( 12\right) \left( 34\right) $ & $\left( 13\right) \left( 24\right) $
& $\left( 14\right) \left( 23\right) $ & $\left( 14\right) \left( 23\right) $
& $\left( 13\right) \left( 24\right) $ & $\left( 14\right) \left( 23\right) $
\\ \hline
$\left( 12\right) \left( 34\right) $ & $\left( 13\right) \left( 24\right) $
& $\left( 14\right) \left( 23\right) $ & $\left( 14\right) \left( 23\right) $
& $\left( 14\right) \left( 23\right) $ & $\left( 13\right) \left( 24\right) $
\\ \hline
$\left( 1234\right) $ & $\left( 1234\right) $ & $\left( 13\right) \left(
24\right) $ & $i$ & $\left( 13\right) $ & $\left( 12\right) \left( 34\right) 
$ \\ \hline
$\left( 1234\right) $ & $\left( 1234\right) $ & $\left( 13\right) \left(
24\right) $ & $i$ & $\left( 12\right) \left( 34\right) $ & $\left( 24\right) 
$ \\ \hline
$\left( 1234\right) $ & $\left( 1432\right) $ & $i$ & $\left( 13\right)
\left( 24\right) $ & $\left( 13\right) $ & $\left( 12\right) \left(
34\right) $ \\ \hline
$\left( 1234\right) $ & $\left( 1432\right) $ & $i$ & $\left( 13\right)
\left( 24\right) $ & $\left( 12\right) \left( 34\right) $ & $\left(
24\right) $ \\ \hline
$\left( 1234\right) $ & $\left( 13\right) $ & $\left( 12\right) \left(
34\right) $ & $\left( 14\right) \left( 23\right) $ & $\left( 12\right)
\left( 34\right) $ & $\left( 24\right) $ \\ \hline
$\left( 1234\right) $ & $\left( 13\right) $ & $\left( 12\right) \left(
34\right) $ & $\left( 14\right) \left( 23\right) $ & $\left( 14\right)
\left( 23\right) $ & $\left( 13\right) $ \\ \hline
$\left( 1234\right) $ & $\left( 12\right) \left( 34\right) $ & $\left(
24\right) $ & $\left( 13\right) $ & $\left( 24\right) $ & $\left( 14\right)
\left( 23\right) $ \\ \hline
$\left( 1234\right) $ & $\left( 12\right) \left( 34\right) $ & $\left(
24\right) $ & $\left( 13\right) $ & $\left( 13\right) $ & $\left( 12\right)
\left( 34\right) $ \\ \hline
\end{tabular}
\end{equation*}

\subsection{\boldmath Figure 15, \guillemotleft $p4gm$\guillemotright}

In the case \guillemotleft $p4gm$\guillemotright , $a$, $b$, $c$, $d$, $u$
and $v$ are as in \guillemotleft $p4$\guillemotright .

Here, $y=cxc^{-1}$.

For $a=i$, $n\leq 2$, the possibilities for $x$ are $i$ and $\left(
12\right) $.

For $a=\left( 12\right) $, $n=2$, the possibilities for $x$ are $i$ and
$\left( 12\right) $.

For $a=\left( 12\right) \left( 34\right) $, $c=\left( 13\right) \left(
24\right) $, the possibilities for $x$ are $\left( 14\right) $ and $\left(
23\right) $.

For $a=\left( 1234\right) $, $c=\left( 1234\right) $, the possibilities for
$x$ are $\left( 13\right) $, $\left( 24\right) $, $\left( 12\right) \left(
34\right) $ and $\left( 14\right) \left( 23\right) $.

For $a=\left( 1234\right) $, $c=\left( 1432\right) $, the possibilities for $x$ are $i$ and $\left( 13\right) \left( 24\right) $.

As \guillemotleft $p4gm$\guillemotright\ is constructed from
\guillemotleft $p4$\guillemotright\ adding reflections, these
reflections can connect permutations:

a) $a=\left( 12\right) $, $c=\left( 34\right) $, $b=\left( 12\right) \left(
34\right) $, $u=i$, $x=\left( 13\right) \left( 24\right) $.

b) $a=\left( 12\right) \left( 34\right) $, $c=\left( 12\right) \left(
34\right) $ , $b=i$, $u=i$, $x=\left( 13\right) \left( 24\right) $.

\begin{equation*}
\begin{tabular}{|c|c|c|c|c|c|}
\hline
$a$ & $c$ & $b$ & $u$ & $x$ & $y$ \\ \hline
$i$ & $i$ & $i$ & $i$ & $i$ & $i$ \\ \hline
$i$ & $i$ & $i$ & $i$ & $\left( 12\right) $ & $\left( 12\right) $ \\ \hline
$\left( 12\right) $ & $\left( 34\right) $ & $\left( 12\right) \left(
34\right) $ & $i$ & $\left( 13\right) \left( 24\right) $ & $\left( 14\right)
\left( 23\right) $ \\ \hline
$\left( 12\right) $ & $\left( 12\right) $ & $i$ & $i$ & $i$ & $i$ \\ \hline
$\left( 12\right) $ & $\left( 12\right) $ & $i$ & $i$ & $\left( 12\right) $
& $\left( 12\right) $ \\ \hline
$\left( 12\right) \left( 34\right) $ & $\left( 12\right) \left( 34\right) $
& $i$ & $i$ & $\left( 13\right) \left( 24\right) $ & $\left( 13\right)
\left( 24\right) $ \\ \hline
$\left( 12\right) \left( 34\right) $ & $\left( 13\right) \left( 24\right) $
& $\left( 14\right) \left( 23\right) $ & $\left( 14\right) \left( 23\right) $
& $\left( 14\right) $ & $\left( 23\right) $ \\ \hline
$\left( 1234\right) $ & $\left( 1234\right) $ & $\left( 13\right) \left(
24\right) $ & $i$ & $\left( 13\right) $ & $\left( 24\right) $ \\ \hline
$\left( 1234\right) $ & $\left( 1234\right) $ & $\left( 13\right) \left(
24\right) $ & $i$ & $\left( 12\right) \left( 34\right) $ & $\left( 14\right)
\left( 23\right) $ \\ \hline
$\left( 1234\right) $ & $\left( 1432\right) $ & $i$ & $\left( 13\right)
\left( 24\right) $ & $i$ & $i$ \\ \hline
$\left( 1234\right) $ & $\left( 1432\right) $ & $i$ & $\left( 13\right)
\left( 24\right) $ & $\left( 13\right) \left( 24\right) $ & $\left(
13\right) \left( 24\right) $ \\ \hline
\end{tabular}
\end{equation*}

\begin{figure}[t]
  \centering
  \includegraphics[width=4.5714in]{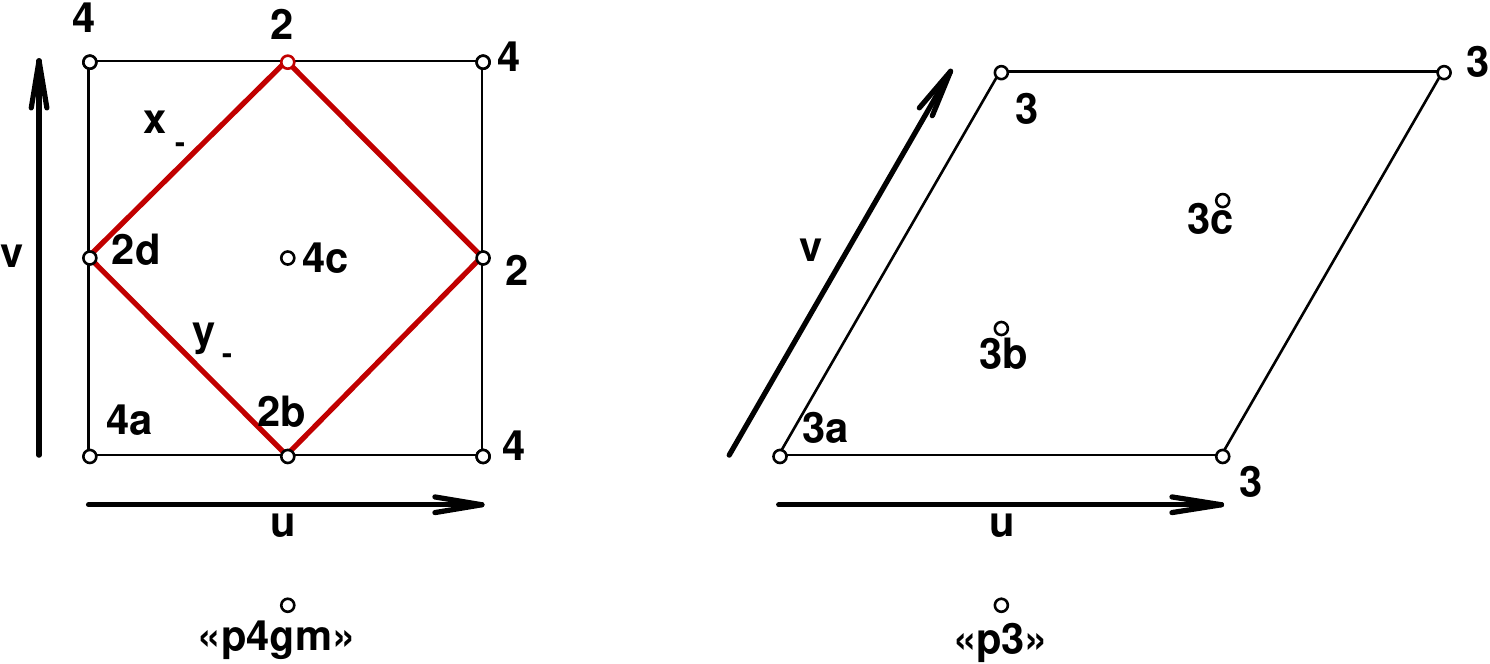}
  \caption{}
  \label{fig:15}
\end{figure}

\subsection{\boldmath Figure 15, \guillemotleft $p3$\guillemotright}

In the case \guillemotleft $p3$\guillemotright , $c=ba^2b$, $u=ba^2$, $v=b^2a $. We shall take always $a$ and $b$ such that $o\left( a\right) \leq
o\left( b\right) ,o\left( c\right) $.

The possibilities for $a$ are $i$ (for $n=1$) and permutations of the type
$\left( 123\right) $ (for $n=3$).

For $a=i$, the possibilities for $b$ are $i$ (for $n=1$) and permutations of
the type $\left( 123\right) $ (for $n=3$).

For $a=\left( 123\right) $, the possibility for $b$ is $\left( 123\right) $.

\begin{equation*}
\begin{tabular}{|c|c|c|c|c|}
\hline
$a$ & $b$ & $c$ & $u$ & $v$ \\ \hline
$i$ & $i$ & $i$ & $i$ & $i$ \\ \hline
$i$ & $\left( 123\right) $ & $\left( 132\right) $ & $\left( 123\right) $ & 
$\left( 132\right) $ \\ \hline
$\left( 123\right) $ & $\left( 123\right) $ & $\left( 123\right) $ & $i$ & $i $ \\ \hline
\end{tabular}
\end{equation*}

\subsection{\boldmath Figure 16, \guillemotleft $p3m1$\guillemotright}

For $a=b=i$, the possibilities for $x$ are $i$ and permutations of the type
$\left( 12\right) $.

For $a=\left( 123\right) $ or $b=\left( 123\right) $, the possibilities for
$x$ are $\left( 12\right) $, $\left( 23\right) $ and $\left( 13\right) $.

\begin{equation*}
\begin{tabular}{|c|c|c|c|c|c|}
\hline
$a$ & $b$ & $c$ & $u$ & $v$ & $x$ \\ \hline
$i$ & $i$ & $i$ & $i$ & $i$ & $i$ \\ \hline
$i$ & $i$ & $i$ & $i$ & $i$ & $\left( 12\right) $ \\ \hline
$i$ & $\left( 123\right) $ & $\left( 132\right) $ & $\left( 123\right) $ &
$\left( 132\right) $ & $i$ \\ \hline
$\left( 123\right) $ & $\left( 123\right) $ & $\left( 123\right) $ & $i$ & 
$i $ & $\left( 12\right) $ \\ \hline
\end{tabular}
\end{equation*}

\begin{figure}[b]
  \centering
  \includegraphics[width=4.6994in]{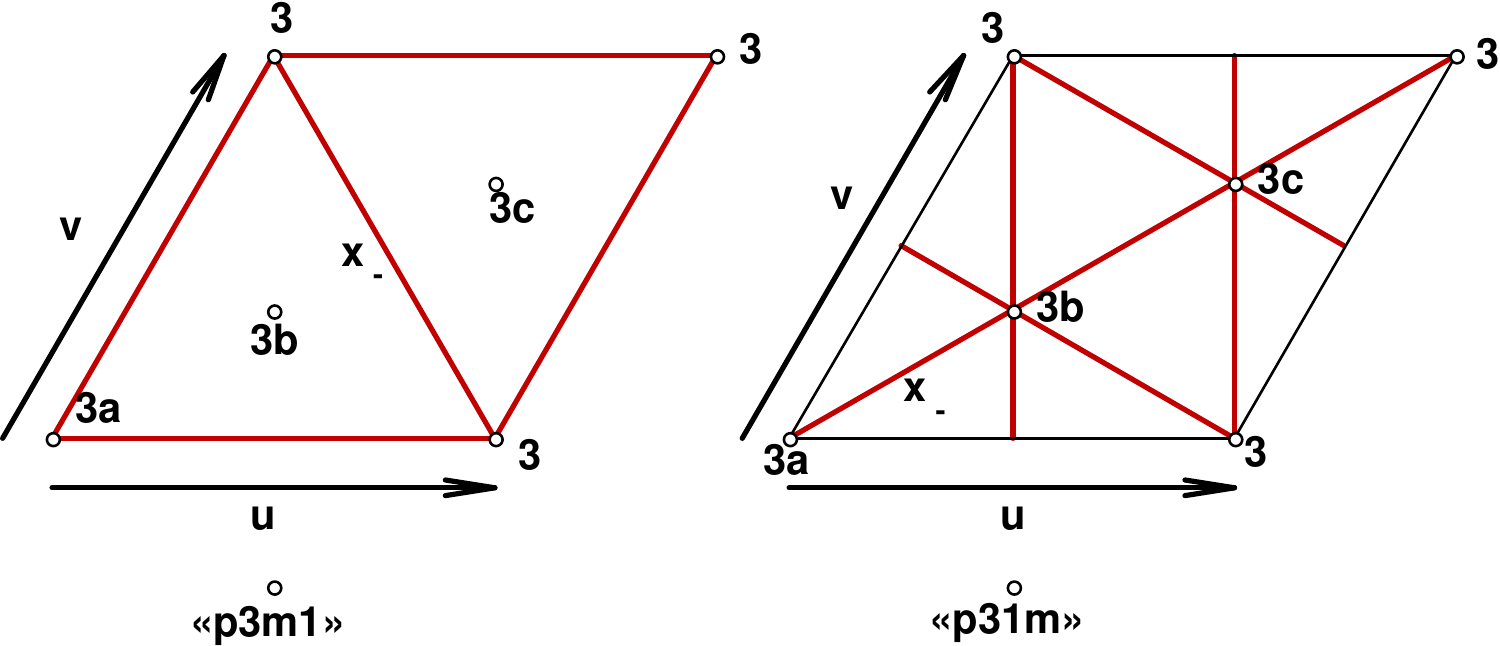}
  \caption{}
  \label{fig:16}
\end{figure}

\subsection{\boldmath Figure 16, \guillemotleft $p31m$\guillemotright}

For $a=b=i$, the possibilities for $x$ are $i$ and permutations of the type 
$\left( 12\right) $.

For $a=\left( 123\right) $ or $b=\left( 123\right) $, the possibilities for
$x$ are $\left( 12\right) $, $\left( 23\right) $ and $\left( 13\right) $.

\begin{equation*}
\begin{tabular}{|c|c|c|c|c|c|}
\hline
$a$ & $b$ & $c$ & $u$ & $v$ & $x$ \\ \hline
$i$ & $i$ & $i$ & $i$ & $i$ & $i$ \\ \hline
$i$ & $i$ & $i$ & $i$ & $i$ & $\left( 12\right) $ \\ \hline
$i$ & $\left( 123\right) $ & $\left( 132\right) $ & $\left( 123\right) $ &
$\left( 132\right) $ & $\left( 12\right) $ \\ \hline
$\left( 123\right) $ & $\left( 123\right) $ & $\left( 123\right) $ & $i$ &
$i $ & $\left( 12\right) $ \\ \hline
\end{tabular}
\end{equation*}

\subsection{\boldmath Figure 17, \guillemotleft $p6$\guillemotright}

In the case \guillemotleft $p6$\guillemotright , $c=ad$, $b=ca$, $u=ca^{-2}$,
$v=c^{-1}a^2$. The permutations $a$ and $d$ are compatible if and only if
they obey the rule $ada=dad$.

For $n=1$, $a=d=i$.

If $o\left( a\right) =o\left( d\right) =2$, the possibilities are the
following:

a) For $n=2$, $a=d=\left( 12\right) $.

b) For $n=3$, $a=\left( 12\right) $, $d=\left( 13\right) $.

c) For $n=6$, $a=\left( 12\right) \left( 34\right) \left( 56\right) $,
$d=\left( 13\right) \left( 25\right) \left( 46\right) $.

If $o\left( a\right) =o\left( d\right) =3$, the possibilities are the
following:

a) For $n=3$, $a=d=\left( 123\right) $.

b) For $n=4$, $a=\left( 123\right) $, $d=\left( 142\right) $.

c) For $n=6$, $a=\left( 123\right) \left( 456\right) $, $d=\left( 124\right)
\left( 356\right) $.

If $o\left( a\right) =o\left( d\right) =6$, the possibilities are the
following:

a) $a=d=\left( 123456\right) $.

b) $a=\left( 123456\right) $, $d=\left( 156423\right) $. There are other two
possibilities for $d$ but they are generated by the action of $a$ on
\textit{this} $d$.

c) $a=\left( 123456\right) $, $d=\left( 163254\right) $. There is other
possibility for $d$ but it is generated by the action of $a$ on \textit{this}
$d$.

\begin{equation*}
\begin{tabular}{|c|c|c|c|c|}
\hline
$a$ & $d$ & $c$ & $b$ & $u$ \\ \hline
$i$ & $i$ & $i$ & $i$ & $i$ \\ \hline
$\left( 12\right) $ & $\left( 12\right) $ & $i$ & $\left( 12\right) $ & $i$
\\ \hline
$\left( 12\right) $ & $\left( 13\right) $ & $\left( 132\right) $ & $\left(
23\right) $ & $\left( 132\right) $ \\ \hline
$\left( 12\right) \left( 34\right) \left( 56\right) $ & $\left( 13\right)
\left( 25\right) \left( 46\right) $ & $\left( 145\right) \left( 263\right) $
& $\left( 16\right) \left( 24\right) \left( 35\right) $ & $\left( 145\right)
\left( 263\right) $ \\ \hline
$\left( 123\right) $ & $\left( 123\right) $ & $\left( 132\right) $ & $i$ &
$i $ \\ \hline
$\left( 123\right) $ & $\left( 142\right) $ & $\left( 143\right) $ & $\left(
12\right) \left( 34\right) $ & $\left( 12\right) \left( 34\right) $ \\ \hline
$\left( 123\right) \left( 456\right) $ & $\left( 124\right) \left(
356\right) $ & $\left( 136\right) \left( 254\right) $ & $\left( 15\right)
\left( 26\right) $ & $\left( 15\right) \left( 26\right) $ \\ \hline
$\left( 123456\right) $ & $\left( 123456\right) $ & $\left( 135\right)
\left( 246\right) $ & $\left( 14\right) \left( 25\right) \left( 36\right) $
& $i$ \\ \hline
$\left( 123456\right) $ & $\left( 156423\right) $ & $\left( 165\right)
\left( 243\right) $ & $\left( 14\right) $ & $\left( 25\right) \left(
36\right) $ \\ \hline
$\left( 123456\right) $ & $\left( 163254\right) $ & $\left( 264\right) $ &
$\left( 16\right) \left( 23\right) \left( 45\right) $ & $\left( 153\right)
\left( 246\right) $ \\ \hline
\end{tabular}
\end{equation*}

\begin{figure}[t]
  \centering
  \includegraphics[width=4.6821in]{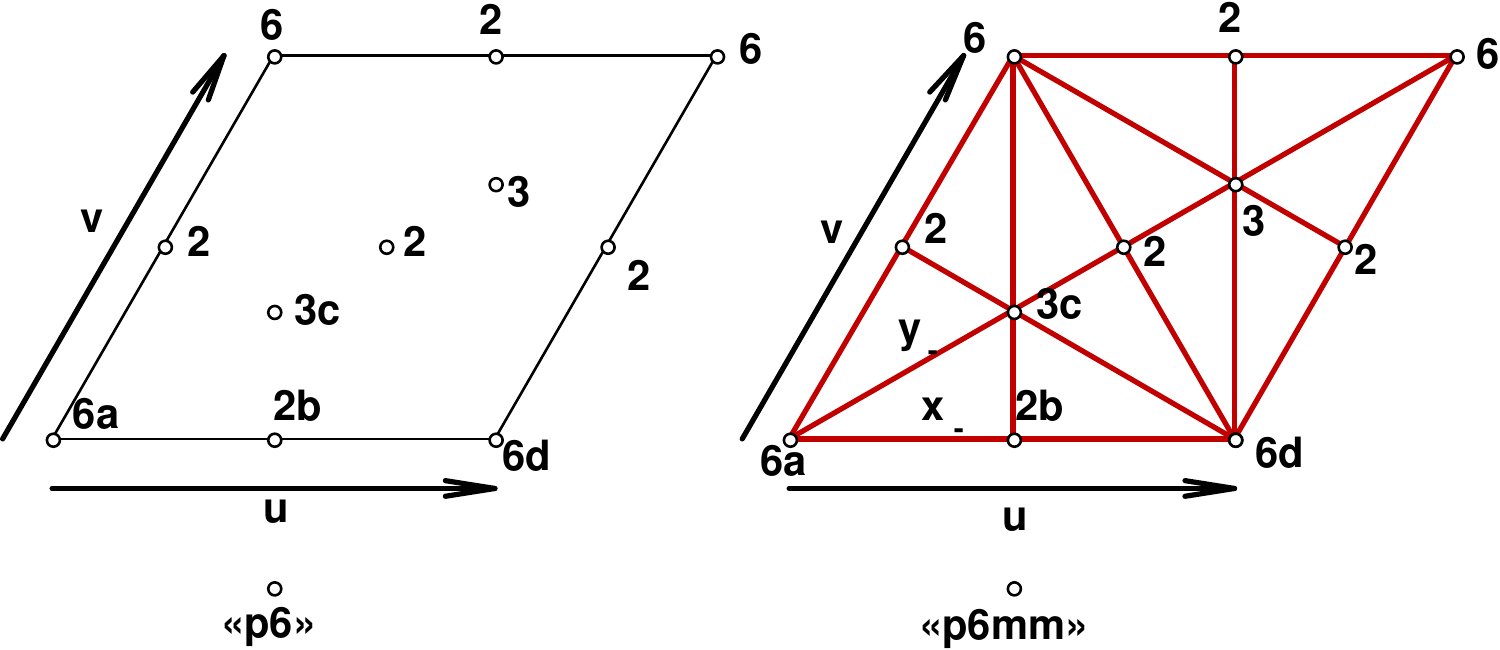}
  \caption{}
  \label{fig:17}
\end{figure}

\subsection{\boldmath Figure 17, \guillemotleft $p6mm$\guillemotright}

In the case \guillemotleft $p6mm$\guillemotright , $a$, $b$, $c$, $d$, $u$
and $v$ are as in \guillemotleft $p6$\guillemotright\ and $y=xa^{-1}$.

As \guillemotleft $p6mm$\guillemotright\ is constructed from
\guillemotleft $p6$\guillemotright\ adding reflections, these reflections can connect
permutations:

a) $a=d=\left( 12\right) \left( 34\right) $. In this case $c=i$, $b=\left(
12\right) \left( 34\right) $, $u=v=i$.

b) $a=\left( 12\right) \left( 45\right) $, $d=\left( 13\right) \left(
46\right) $. In this case $c=\left( 132\right) \left( 465\right) $,
$b=\left( 23\right) \left( 56\right) $, $u=\left( 132\right) \left(
465\right) $, $v=\left( 123\right) \left( 456\right) $.

c) $a=d=\left( 123\right) \left( 456\right) $. In this case $c=\left(
132\right) \left( 465\right) $, $b=i$, $u=v=i$.

\begin{equation*}
\begin{tabular}{|c|c|c|c|}
\hline
$a$ & $d$ & $x$ & $y$ \\ \hline
$i$ & $i$ & $i$ & $i$ \\ \hline
$i$ & $i$ & $\left( 12\right) $ & $\left( 12\right) $ \\ \hline
$\left( 12\right) $ & $\left( 12\right) $ & $i$ & $\left( 12\right) $ \\ 
\hline
$\left( 12\right) $ & $\left( 12\right) $ & $\left( 12\right) $ & $i$ \\ 
\hline
$\left( 12\right) \left( 34\right) $ & $\left( 12\right) \left( 34\right) $
& $\left( 13\right) \left( 24\right) $ & $\left( 14\right) \left( 23\right) $
\\ \hline
$\left( 12\right) $ & $\left( 13\right) $ & $i$ & $\left( 12\right) $ \\ 
\hline
$\left( 12\right) \left( 45\right) $ & $\left( 13\right) \left( 46\right) $
& $\left( 14\right) \left( 25\right) \left( 36\right) $ & $\left( 15\right)
\left( 24\right) \left( 36\right) $ \\ \hline
$\left( 12\right) \left( 34\right) \left( 56\right) $ & $\left( 13\right)
\left( 25\right) \left( 46\right) $ & $i$ & $\left( 12\right) \left(
34\right) \left( 56\right) $ \\ \hline
$\left( 12\right) \left( 34\right) \left( 56\right) $ & $\left( 13\right)
\left( 25\right) \left( 46\right) $ & $\left( 12\right) \left( 35\right)
\left( 46\right) $ & $\left( 36\right) \left( 45\right) $ \\ \hline
$\left( 12\right) \left( 34\right) \left( 56\right) $ & $\left( 13\right)
\left( 25\right) \left( 46\right) $ & $\left( 16\right) \left( 25\right)
\left( 34\right) $ & $\left( 15\right) \left( 26\right) $ \\ \hline
$\left( 12\right) \left( 34\right) \left( 56\right) $ & $\left( 13\right)
\left( 25\right) \left( 46\right) $ & $\left( 13\right) \left( 24\right)
\left( 56\right) $ & $\left( 14\right) \left( 23\right) $ \\ \hline
$\left( 123\right) $ & $\left( 123\right) $ & $\left( 12\right) $ & $\left(
13\right) $ \\ \hline
$\left( 123\right) \left( 456\right) $ & $\left( 123\right) \left(
456\right) $ & $\left( 14\right) \left( 26\right) \left( 35\right) $ & $\left( 15\right) \left( 24\right) \left( 36\right) $ \\ \hline
$\left( 123\right) $ & $\left( 142\right) $ & $\left( 12\right) $ & $\left(
13\right) $ \\ \hline
$\left( 123\right) \left( 456\right) $ & $\left( 124\right) \left(
356\right) $ & $\left( 12\right) \left( 56\right) $ & $\left( 13\right)
\left( 45\right) $ \\ \hline
$\left( 123\right) \left( 456\right) $ & $\left( 124\right) \left(
356\right) $ & $\left( 16\right) \left( 25\right) \left( 34\right) $ & $\left( 14\right) \left( 26\right) \left( 35\right) $ \\ \hline
$\left( 123456\right) $ & $\left( 123456\right) $ & $\left( 14\right) \left(
23\right) \left( 56\right) $ & $\left( 15\right) \left( 24\right) $ \\ \hline
$\left( 123456\right) $ & $\left( 123456\right) $ & $\left( 26\right) \left(
35\right) $ & $\left( 12\right) \left( 36\right) \left( 45\right) $ \\ \hline
$\left( 123456\right) $ & $\left( 156423\right) $ & $\left( 14\right) \left(
23\right) \left( 56\right) $ & $\left( 15\right) \left( 24\right) $ \\ \hline
$\left( 123456\right) $ & $\left( 156423\right) $ & $\left( 26\right) \left(
35\right) $ & $\left( 12\right) \left( 36\right) \left( 45\right) $ \\ \hline
$\left( 123456\right) $ & $\left( 163254\right) $ & $\left( 14\right) \left(
23\right) \left( 56\right) $ & $\left( 15\right) \left( 24\right) $ \\ \hline
\end{tabular}
\end{equation*}

\section{Plane symmetries, permutation groups and puzzle solutions}

\subsection{Board symmetries}

Consider a board with at least three faces with a common vertex. From now
on $V$ denotes the set of the board vertices, $E$ denotes the set of the
board edges and $F$ denotes the set of the board faces (polygons).

The group of the board symmetries, $\Omega $, called the board group, is the
set of all isometries $\omega $ of $\mathbb{R}^2$, $\omega \equiv \left( u,\eta
\right) $, that send vertices to vertices, which implies that they send
edges to edges, faces to faces. Every symmetry $\omega \in \Omega $\ induces
three bijections, that we shall also denote $\omega $, whenever there is no
confusion possible: $\omega :V\rightarrow V$, $\omega :E\rightarrow E$ and
$\omega :F\rightarrow F$. Denote also $\Omega \equiv \left\{ \omega
:V\rightarrow V\right\} \equiv \left\{ \omega :E\rightarrow E\right\} \equiv
\left\{ \omega :F\rightarrow F\right\} $, the three sets of these functions.
One can say that each one of these three sets $\Omega $\ is the set of the
board symmetries. Notice that not all one-to-one functions $F\longrightarrow
F$, $E\longrightarrow E$, $V\longrightarrow V$ are in $\Omega $. With the
composition of functions each one of these three sets $\Omega $ forms a
group that is isomorphic to the group of the board symmetries. If $\omega
_1,\omega _2\in \Omega $, we shall denote $\omega _1\omega _2\equiv \omega
_1\circ \omega _2$.

When no confusion is possible, $\omega \in \Omega $ represents also the
group isomorphism $\omega :\Omega \rightarrow \Omega $, $\omega \left(
\omega _1\right) =\omega \omega _1\omega ^{-1}$, for every $\omega _1\in
\Omega $. Note that $\omega _1$ and $\omega \left( \omega _1\right) $\ have
the same order.

If $\Omega _{1}$ is a subgroup of $\Omega $, then $\Omega _{1}$ acts
naturally on the face set, $F$: for $\omega \in \Omega _{1}$ and $\varphi
\in F$, one defines the action $\omega \varphi =\omega \left( \varphi
\right) $.

\subsection{Puzzle solutions}

Consider a puzzle with numbers $1,2,\ldots ,n$ drawn on the plates. From now
on $P$ denotes the set of its plates which have numbers drawn, and call it
the plate set. If no confusion is possible, $P$ will also denote the puzzle
itself. However, note that we can not separate the plates from the board:
the puzzle is the plates \textit{and} the board.

As before $E$ denotes the set of the board edges and $F$ denotes the set of
the board faces. A solution of the puzzle defines a function $\varepsilon
:E\rightarrow \left\{ 1,2,\ldots ,n\right\} $. Denote $\mathcal{E}$ the set
of these functions. One can say that $\mathcal{E}$\ is the set of the puzzle
solutions.

We shall also consider the group $S_n\times \Omega $. If $\left( a_1,\omega
_1\right) ,\left( a_2,\omega _2\right) \in S_n\times \Omega $, one defines
the product $\left( a_1,\omega _1\right) \left( a_2,\omega _2\right) =\left(
a_1a_2,\omega _1\omega _2\right) $.

\subsection{The plate group}

Some $S_n$\ subgroups act naturally on $P$. Let $\pi \in P$ and $a\in S_n$.
Assume that $m_1,m_2,m_3,\ldots $ are drawn on $\pi $, by this order. Then
$a\pi $ is a plate where the numbers $a\left( m_1\right) =n_1,a\left(
m_2\right) =n_2,a\left( m_3\right) =n_3,\ldots $ are drawn replacing
$m_1,m_2,m_3,\ldots $ (see Figure 18).

\begin{figure}
  \centering
  \includegraphics[width=2.0617in]{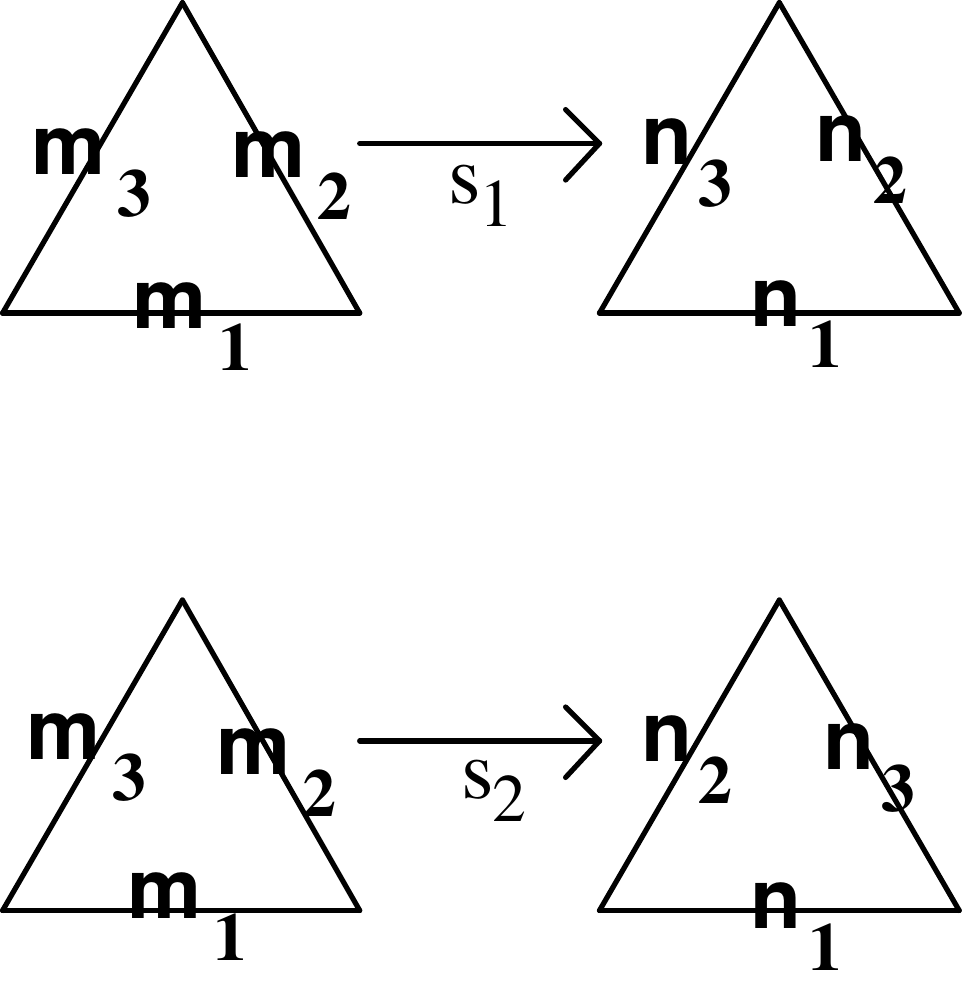}
  \caption{}
  \label{fig:18}
\end{figure}

Let $s\in S_n^{\pm }$ and $\pi \in P$. If $s\equiv s_1=\left( 1,a\right)
\equiv a$, then $s\pi =a\pi $. If $s\equiv s_2=\left( -1,a\right) \equiv
a^{-}$, then $s\pi $ is a reflection of $a\pi $. In this last case, if the
numbers $m_1,m_2,m_3,\ldots $ are drawn on $\pi $, by this order, then $s\pi 
$ is a plate where the numbers $\ldots ,a\left( m_3\right) =n_3,a\left(
m_2\right) =n_2,a\left( m_1\right) =n_1$ are drawn by this order (see Figure
18).

The plate group, $G_P$, is the greatest subgroup of $S_n^{\pm }$ (if $\Omega
^{-}\neq \emptyset $) that acts on $P$, or the greatest subgroup of $S_n$
(if $\Omega ^{-}=\emptyset $) that acts on $P$. For $\Omega ^{-}\neq
\emptyset $, if $s\in S_n^{\pm }$ and $s\pi \in P$, for every $\pi \in P$,
then $s\in G_P$. For $\Omega ^{-}=\emptyset $, if $s\in S_n$ and $s\pi \in P 
$, for every $\pi \in P$, then $s\in G_P$.

\subsection{The solution group}

Let $\varepsilon :E\rightarrow \left\{ 1,2,\ldots ,n\right\} $ be a solution
of the puzzle. The group of this solution, $G_{\varepsilon _{}}$, is a
subgroup of $S_n\times \Omega $; $\left( a,\omega \right) \in G_{\varepsilon
_{}}$ if and only if 
\begin{equation*}
a\circ \varepsilon =\varepsilon \circ \omega \text{.}
\end{equation*}

Denote $\Omega _{\varepsilon _{}}$ the following subgroup of $\Omega $:
$\omega \in \Omega _{\varepsilon _{}}$ if and only if there exists $a\in S_n$
such that $\left( a,\omega \right) \in G_{\varepsilon _{}}$. Notice that if
$\omega \in \Omega _{\varepsilon _{}}$ there exists only one $a\in S_n$ such
that $\left( a,\omega \right) \in G_{\varepsilon _{}}$. From this one
concludes that $\omega \mapsto \left( a,\omega \right) $ defines an
isomorphism between $\Omega _{\varepsilon _{}}$ and $G_{\varepsilon _{}}$
and that $\left( \det \eta ,a\right) \in G_P$. This defines $g_{\varepsilon
_{}}:\Omega _{\varepsilon _{}}\rightarrow G_P$, $\ g_{\varepsilon _{}}\left(
\omega \right) =\left( \det \eta ,a\right) $, which is an homomorphism of
groups.

For a lot of puzzles $\left( \det \eta ,a\right) $ defines completely
$\omega $. It is the case of all puzzles considered in this article. Hence,
when $\left( \det \eta ,a\right) $ defines completely $\omega $,
$g_{\varepsilon _{}}$ establishes an isomorphism between $\Omega
_{\varepsilon _{}}$ and $g_{\varepsilon _{}}\left( \Omega _{\varepsilon
_{}}\right) \subset G_P$. Denote $G_{P_{\varepsilon _{}}}\equiv
g_{\varepsilon _{}}\left( \Omega _{\varepsilon _{}}\right) $. Finally,
$G_{\varepsilon _{}}$ and $G_{P_{\varepsilon _{}}}$ are isomorphic. We can
identify $\left( a,\omega \right) $ with $\left( \det \eta ,a\right) $, and
$G_{\varepsilon _{}}$ with the subgroup $G_{P_{\varepsilon _{}}}$ of $G_P$.

\subsection{Equivalent solutions}

Let $\varepsilon _1,\varepsilon _2:E\rightarrow \left\{ 1,2,\ldots
,n\right\} $ be solutions of the puzzle. One says that these solutions are
equivalent, $\varepsilon _1\approx \varepsilon _2$, if there are $\omega \in
\Omega $, $\omega \equiv \left( u,\eta \right) $,and $a\in S_n$, such that 
\begin{equation*}
a\circ \varepsilon _1=\varepsilon _2\circ \omega \text{.}
\end{equation*}

Notice that $\left( \det \eta ,a\right) \in G_P$.

If $a=i$ and $\det \eta =1$, what distinguishes the solutions $\varepsilon
_1 $ and $\varepsilon _2$ is only a symmetry in $\Omega ^{+}$. In this case 
\begin{equation*}
\varepsilon _1=\varepsilon _2\circ \omega
\end{equation*}
expresses another equivalence relation, $\varepsilon _1\sim \varepsilon _2$.
When we make a puzzle, in practice, we do not recognize the difference
between $\varepsilon _1$ and $\varepsilon _2$. We shall say that they
represent the same \textit{natural solution}, an equivalence class of the
relation $\sim $.

Let $\varepsilon ,\varepsilon _{1},\varepsilon _{2}\in \mathcal{E}$. As
$\varepsilon _{1}\sim \varepsilon _{2}$ and $\varepsilon _{1}\approx
\varepsilon $, implies $\varepsilon _{2}\approx \varepsilon $, one can say
that the natural solution represented by $\varepsilon _{1}$ is equivalent to 
$\varepsilon $.

This equation involving $\varepsilon _1$ and $\varepsilon _2$ defines an
equivalence relation, and a natural solution is an equivalence class of this
relation. Notice that if $\varepsilon _1=\varepsilon _2$, then $\omega _1$
is the identity.

For $\varepsilon \in \mathcal{E}$, represent by $\left[ \varepsilon \right] $\ the set of natural solutions equivalent to $\varepsilon $.

Now, if $\Omega ^{-}\neq \emptyset $, choose $\omega _{-}\in \Omega ^{-}$.
For $\varepsilon \in \mathcal{E}$ and $s=\left( \delta ,a\right) \in G_P$,
denote $\varepsilon _s=a\circ \varepsilon \circ \omega $, where $\omega $ is
the identity if $\delta =1$ and $\omega =\omega _{-}$ if $\delta =-1$. The
set $\left\{ \varepsilon _s:s\in G_P\right\} $ includes representatives of
all natural solutions equivalent to $\varepsilon $. Then 
\begin{equation*}
\left| \left[ \varepsilon \right] \right| =\frac{\left| G_P\right| }{\left|
G_{P_{\varepsilon _{}}}\right| }\text{.}
\end{equation*}

The cardinal of all the natural solutions is then given by
\begin{equation*}
\sum_{\left[ \varepsilon \right] }\frac{\left| G_P\right| }{\left|
G_{P_{\varepsilon _{}}}\right| }\text{,}
\end{equation*}
where the sum is extended to all different equivalence classes $\left[
\varepsilon \right] $.

\subsection{Equivalent puzzles}

Consider two puzzles and their plate sets, $P_1$ and $P_2$. Let $\delta =1$
if they have the same board, but $\Omega ^{-}=\emptyset $; $\delta =-1$ if
they have different boards, but symmetric by reflection; and $\delta =\pm 1$
if they have the same board and $\Omega ^{-}\neq \emptyset $. One says that
the puzzles are equivalent if there exists $s=\left( \delta ,a\right) \in
S_n^{\pm }$ such that the function $\pi \longmapsto s\pi $ is one-to-one
between $P_1$ and $P_2$. We denote $P_2$ as $sP_1$.

As an example take the puzzles in Figure 20 (c) and (d). Let $a=\left(
13\right) \left( 67\right) $. If $s=\left( -1,a\right) $ acts on the plates
of one of these puzzles, it generates the plates of the other. If $P$ is the
plate set of one of these puzzles, then $P\neq sP$.

\section{Examples}

In this section we use group theory in order to find puzzles, for a given
board, such as, for example, maximal puzzles. These are important examples,
but others could be given.

To avoid ambiguities, in the puzzles we give in the following, all the edges
have numbers.

Some of the examples we give in this section can easily be studied directly.
All the results we present are obtained in this way. Some puzzles have only
one natural solution. Others can have millions of natural solutions that can
only be calculated with a computer.

\subsection{\boldmath Puzzles \guillemotleft $p4$\guillemotright} \ 

\vspace{-10pt}
\begin{figure}[H]
  \centering
  \includegraphics[width=5.78in]{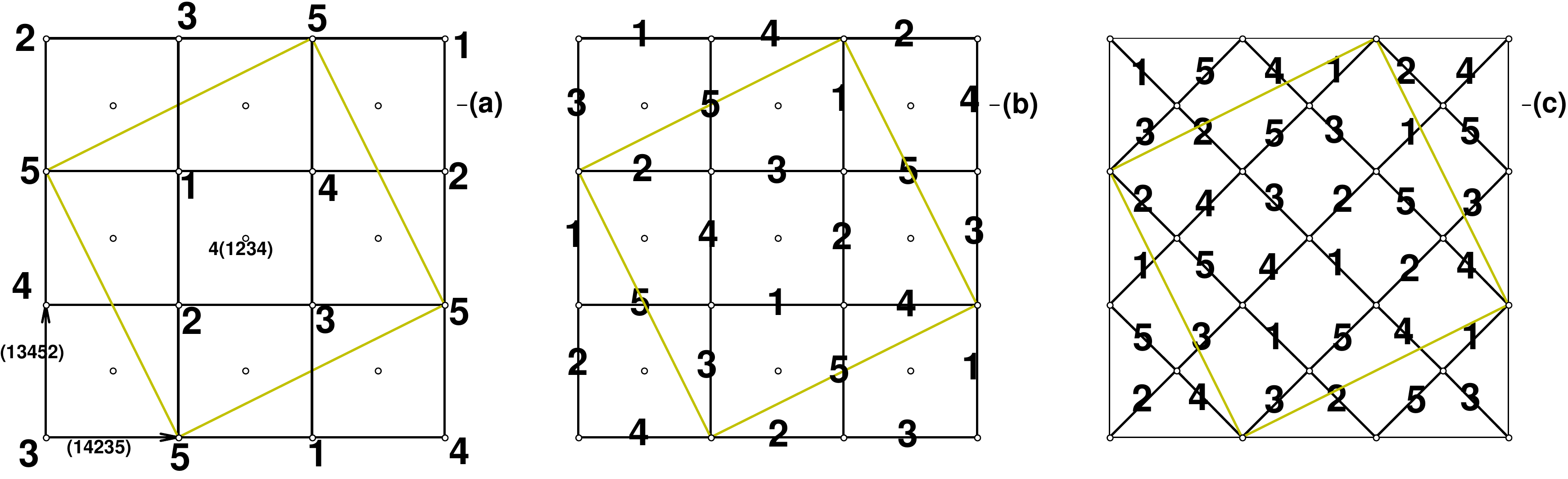}
  \vspace{-30pt}
  \caption{}
  \label{fig:19}
\end{figure}

\subsubsection{Figure 19 (a)}

Take a board with five square faces (Figure 19 (a)). Obviously, it is not
symmetric by reflection. The group of this board has $20$ elements.

The group of permutations in $S_5$ associated with translations of this
pattern are generated by $u=\left( 14235\right) $ and $v=\left( 13452\right) 
$. The group of this pattern is generated by 
\begin{equation*}
\begin{tabular}{|c|c|c|c|c|}
\hline
$a$ & $c$ & $b$ & $u$ & $v$ \\ \hline
$\left( 1234\right) $ & $\left( 1325\right) $ & $\left( 15\right) \left(
34\right) $ & $\left( 14235\right) $ & $\left( 13452\right) $ \\ \hline
\end{tabular}
\text{.}
\end{equation*}

\subsubsection{Figure 19 (b)}

In this figure it is represented a maximal solution of a puzzle with plates
$\left[ 1234\right] $, $\left[ 1325\right] $, $\left[ 1453\right] $, $\left[
1542\right] $, $\left[ 2435\right] $. It is not very interesting because, in
fact, this natural solution is unique. Its group has $20$ elements.

\subsubsection{Figure 19 (c)}

One can construct another board with ten square faces with vertices in the
ten rotations centers of the grid (Figure 19 (c)). The group of this board
has $40$ elements.

In Figure 19 (c), is represented a maximal solution of a puzzle with plates
$\left[ 1415\right] $, $\left[ 1213\right] $, $\left[ 2125\right] $, $\left[
2324\right] $, $\left[ 3134\right] $, $\left[ 3235\right] $, $\left[
4142\right] $, $\left[ 4345\right] $, $\left[ 5153\right] $, $\left[
5254\right] $. Its group has $20$ elements.

\subsection{\boldmath Puzzles \guillemotleft $p4mm$\guillemotright}

\subsubsection{Figure 20 (a)}

The board with eight faces represented in Figure 20 (a) is symmetric by
reflection. The group of this board has $64$ elements.

The group of permutations in $S_8$ associated with translations of the
pattern in the point lattice is generated by $u=\left( 1835\right) \left(
2647\right) $ and $v=\left( 1637\right) \left( 2845\right) $. The group of
this pattern is generated by 
\begin{equation*}
\begin{tabular}{|c|c|c|}
\hline
$a$ & $c$ & $b$ \\ \hline
$\left( 1234\right) \left( 67\right) $ & $\left( 1265\right) \left(
3748\right) $ & $\left( 15\right) \left( 27\right) \left( 38\right) \left(
46\right) $ \\ \hline
\end{tabular}
\end{equation*}
\begin{equation*}
\begin{tabular}{|c|c|c|}
\hline
$u$ & $x$ & $y$ \\ \hline
$\left( 1835\right) \left( 2647\right) $ & $\left( 13\right) \left(
67\right) $ & $\left( 12\right) \left( 34\right) $ \\ \hline
\end{tabular}
\end{equation*}

There are no puzzles with solutions that are ``compatible'' with 
$\Omega^{+} $.

One can construct another board with sixteen square faces with vertices in
the sixteen rotations centers of the grid. The group of this board has $128$
elements.

There are three puzzles with solutions that are ``compatible'' with $\Omega
^{+}$. In Figure 20 (b), (c) and (d) we represent their unique maximal
solutions. Puzzles (c) and (d) are equivalent.

\begin{figure}[t]
  \centering
  \includegraphics[width=4.7703in]{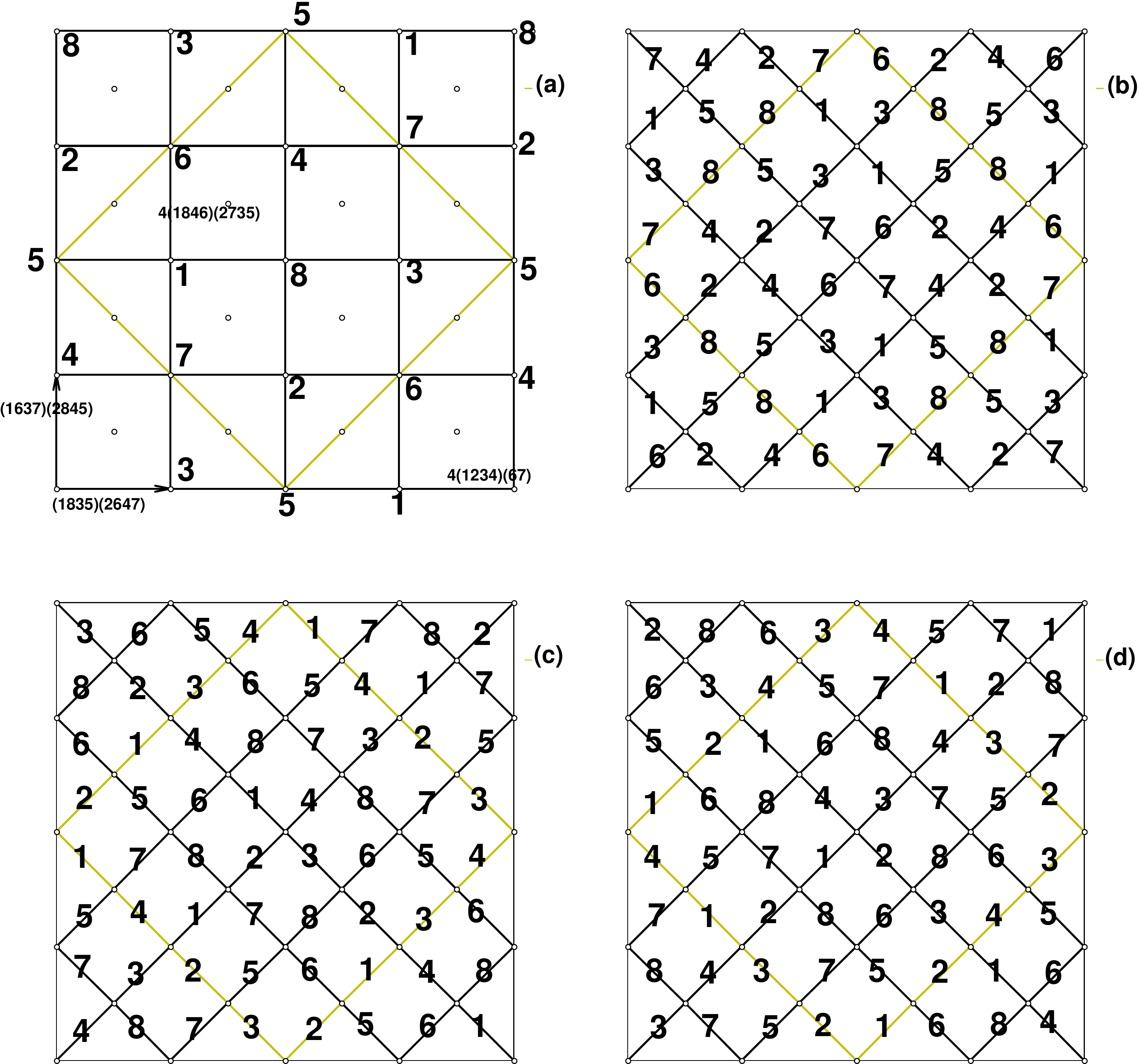}
  \caption{}
  \label{fig:20}
\end{figure}

\subsubsection{Figure 20 (b)}

In this figure it is represented a maximal solution of a puzzle with plates
$\left[ 1358\right] $, $\left[ 1367\right] $, $\left[ 1376\right] $, $\left[
1385\right] $, $\left[ 1542\right] $, $\left[ 1583\right] $, $\left[
1673\right] $, $\left[ 1763\right] $, $\left[ 1853\right] $, $\left[
2458\right] $, $\left[ 2467\right] $, $\left[ 2476\right] $, $\left[
2485\right] $, $\left[ 2674\right] $, $\left[ 2764\right] $, $\left[
2854\right] $. Its group has $64$ elements.

\subsubsection{Figure 20 (c)}

In this figure it is represented a maximal solution of a puzzle with plates
$\left[ 1257\right] $, $\left[ 1275\right] $, $\left[ 1286\right] $, $\left[
1465\right] $, $\left[ 1478\right] $, $\left[ 1564\right] $, $\left[
1682\right] $, $\left[ 1752\right] $,$\left[ 1874\right] $, $\left[
2356\right] $, $\left[ 2387\right] $, $\left[ 2783\right] $, $\left[
2653\right] $, $\left[ 3457\right] $, $\left[ 3486\right] $, $\left[
3754\right] $. Its group has $32$ elements.

\subsubsection{Figure 20 (d)}

The other puzzle in Figure 20 (d) is equivalent to the puzzle in Figure 20
(c). For example, if $\left( 13\right) \left( 67\right) _{-}$ acts on the
plates of one of these puzzles, it generates the plates of the other.

\subsection{\boldmath Puzzles \guillemotleft $p6$\guillemotright}

\subsubsection{Figure 21}

\begin{figure}[b]
  \centering
  \includegraphics[width=5.6368in]{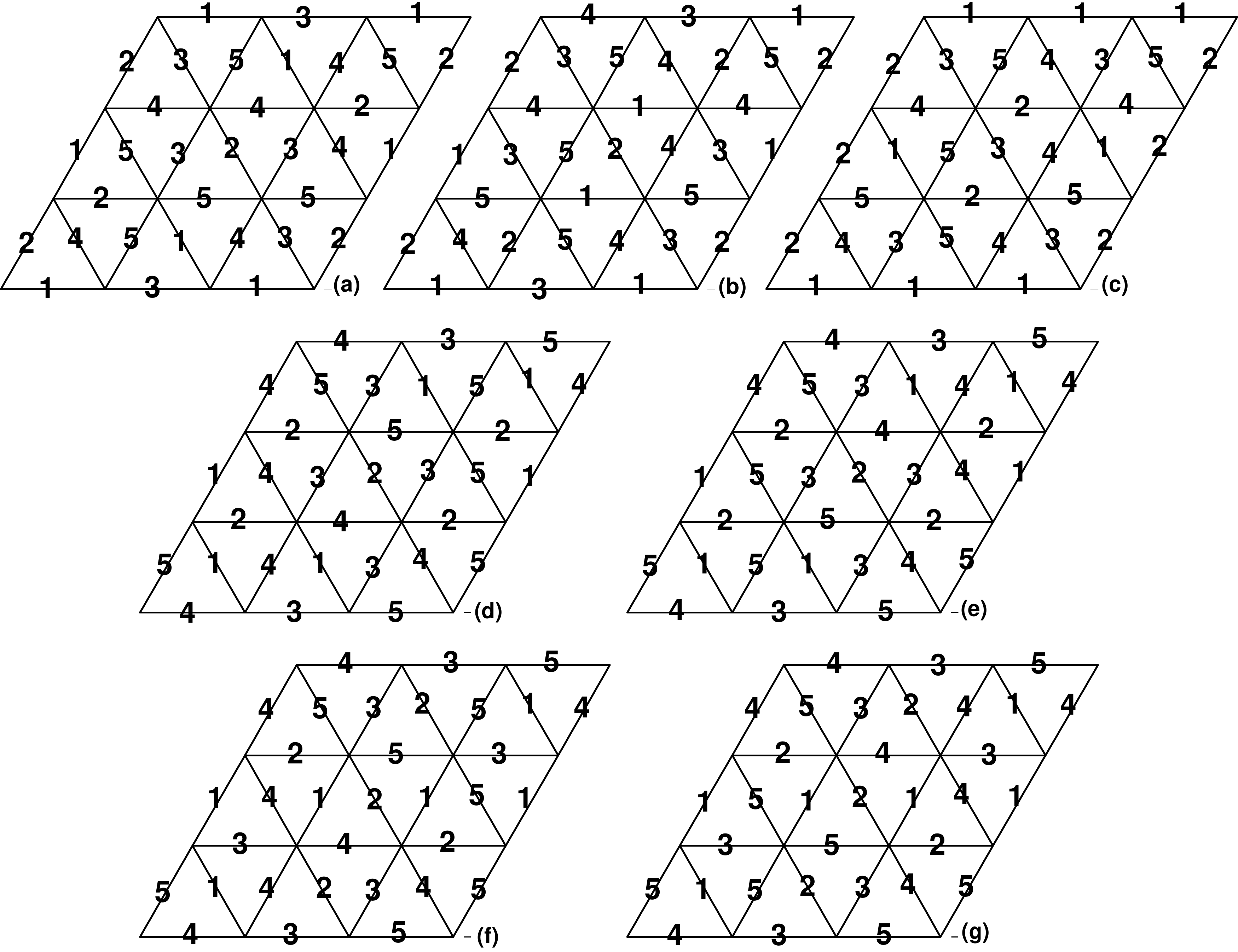}
  \caption{}
  \label{fig:21}
\end{figure}

The non connected group is generated by 
\begin{equation*}
\begin{tabular}{|c|c|c|c|c|}
\hline
$a$ & $d$ & $c$ & $b$ & $u$ \\ \hline
$\left( 123\right) \left( 45\right) $ & $\left( 123\right) \left( 45\right) $
& $\left( 132\right) $ & $\left( 45\right) $ & $i$ \\ \hline
\end{tabular}
\end{equation*}

The plates are: $\left[ 124\right] $, $\left[ 125\right] $, $\left[
134\right] $, $\left[ 135\right] $, $\left[ 142\right] $, $\left[ 143\right] 
$, $\left[ 145\right] $, $\left[ 152\right] $, $\left[ 153\right] $, $\left[
154\right] $, $\left[ 234\right] $, $\left[ 235\right] $, $\left[ 243\right] 
$, $\left[ 245\right] $, $\left[ 253\right] $, $\left[ 254\right] $, $\left[
345\right] $, $\left[ 354\right] $.

In this puzzle $\left| G_P\right| =24$. Figure 21 shows natural solutions
that represent the $7$ equivalence classes a rotation of order $6$. Each one
has $4$ natural solutions and a group of order $6$.

\subsubsection{Figure 22}

Take a board with $14$ triangular equilateral faces (see Figure 22 (a)).
Obviously, it is not symmetric by reflection. The group of this board has
$42$ elements.

The group of permutations in $S_7$ associated with translations of this
pattern are generated by $\left( 1325647\right) $ and $\left( 1572436\right) 
$.

\begin{figure}[b]
  \centering
  \includegraphics[width=5.7372in]{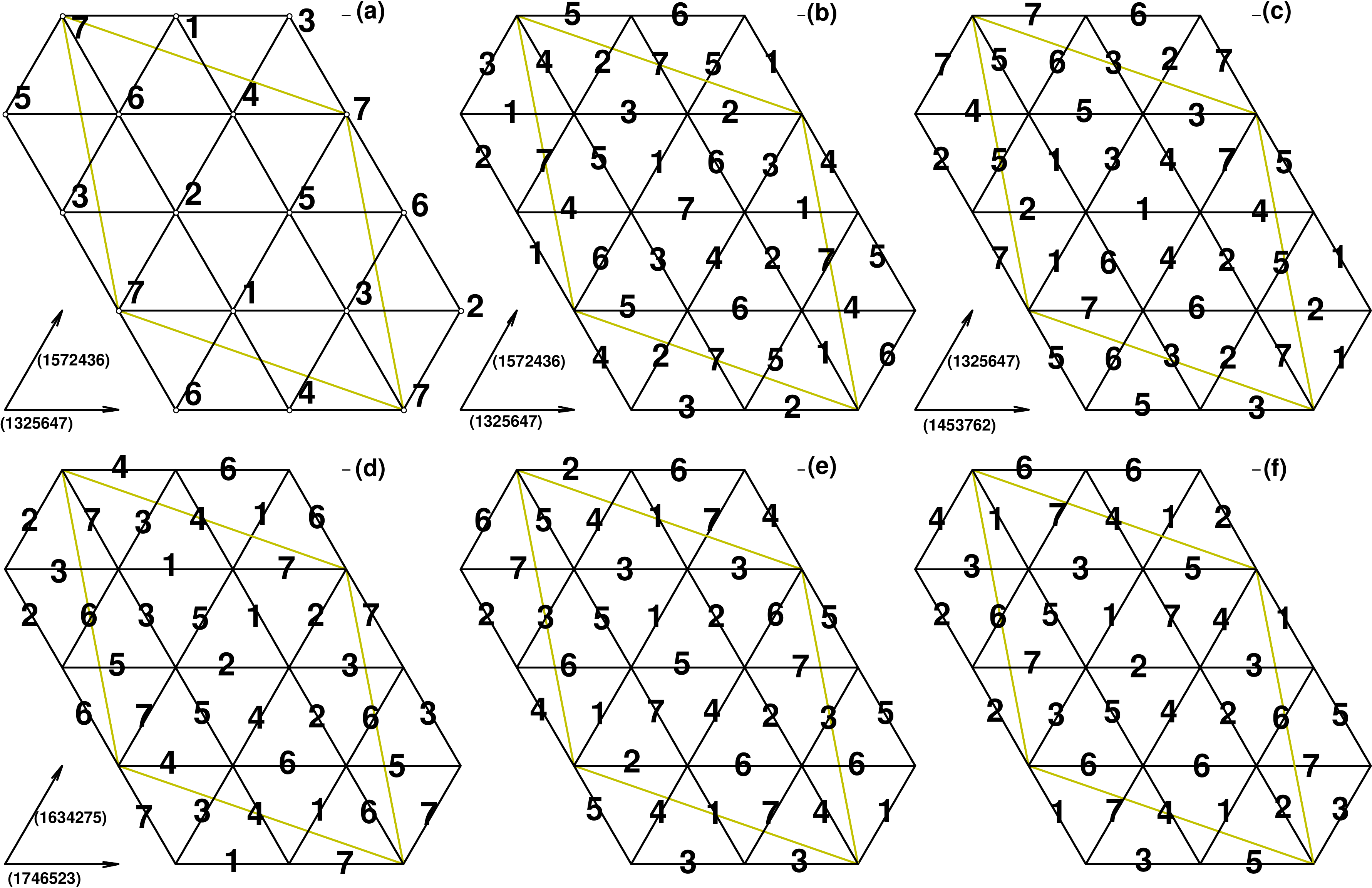}
  \vspace{-12pt}
  \caption{}
  \label{fig:22}
\end{figure}

The numbers over the point lattice generate a pattern. Associated with the
rotations of order $6$ (angle of $\pi /3$ in the direct sense) centered in
the vertices of the grid (the points of the lattice) are permutations like
$\left( 123456\right) $ and all the others generated by the translations and
this one. Associated with the rotations of order $3$ (angle of $2\pi /3$ in
the direct sense) centered in the triangles are permutations like $\left(
135\right) \left( 246\right) $ and all the others generated by the
translations and this one.

The group of the pattern is generated by 
\begin{equation*}
\begin{tabular}{|c|c|c|c|c|}
\hline
$a$ & $d$ & $c$ & $b$ & $u$ \\ \hline
$\left( 123456\right) $ & $\left( 276435\right) $ & $\left( 127\right)
\left( 365\right) $ & $\left( 17\right) \left( 26\right) \left( 34\right) $
& $\left( 1325647\right) $ \\ \hline
\end{tabular}
\end{equation*}

We naturally construct triangular plates based upon these permutations of
order $3$: $\left[ 127\right] $, $\left[ 135\right] $, $\left[ 143\right] $, 
$\left[ 152\right] $, $\left[ 164\right] $, $\left[ 176\right] $, $\left[
237\right] $, $\left[ 246\right] $, $\left[ 254\right] $, $\left[ 263\right] 
$, $\left[ 365\right] $, $\left[ 347\right] $, $\left[ 457\right] $, $\left[
567\right] $. In Figure 22 (b) is represented the maximal natural solution
of this puzzle. Its group has $42$ elements.

We fix now plate $\left[ 246\right] $ in the board and look for all
solutions that have associated the permutation $\left( 135\right) \left(
246\right) $ to the rotation of order $3$ in the middle of this triangle.
There are $5$ such solutions: (b)--(f).

The solutions (c) and (d) are equivalent. One can obtain (d) from (c) by
operating with the permutation $\left( 123456\right) $ over the plates.
Their group has no permutations of order $6$, and it has $21$ elements.

The solutions (e) and (f) are equivalent. One can obtain (f) from (e) by
operating with the permutation $\left( 14\right) \left( 25\right) \left(
36\right) $ over the plates. Their group has no permutations of order $6$,
and it has $3$ elements.

\subsubsection{Figure 23}

The group is generated by 
\begin{equation*}
\begin{tabular}{|c|c|c|c|c|}
\hline
$a$ & $d$ & $c$ & $b$ & $u$ \\ \hline
$\left( 123456\right) $ & $\left( 156423\right) $ & $\left( 165\right)
\left( 243\right) $ & $\left( 14\right) $ & $\left( 25\right) \left(
36\right) $ \\ \hline
\end{tabular}
\end{equation*}

In this figure we have a board of eight equilateral triangles as faces. The
plates are $\left[ 132\right] $, $\left[ 135\right] $, $\left[ 162\right] $, 
$\left[ 165\right] $, $\left[ 243\right] $, $\left[ 246\right] $, $\left[
354\right] $, $\left[ 465\right] $. There are three equivalence classes with 
$2$ (a) $6$ (b) and $6$ (c) natural solutions. Figure 23 includes
representatives of the three equivalence classes.

\begin{figure}[H]
  \centering
  \includegraphics[width=5.271in]{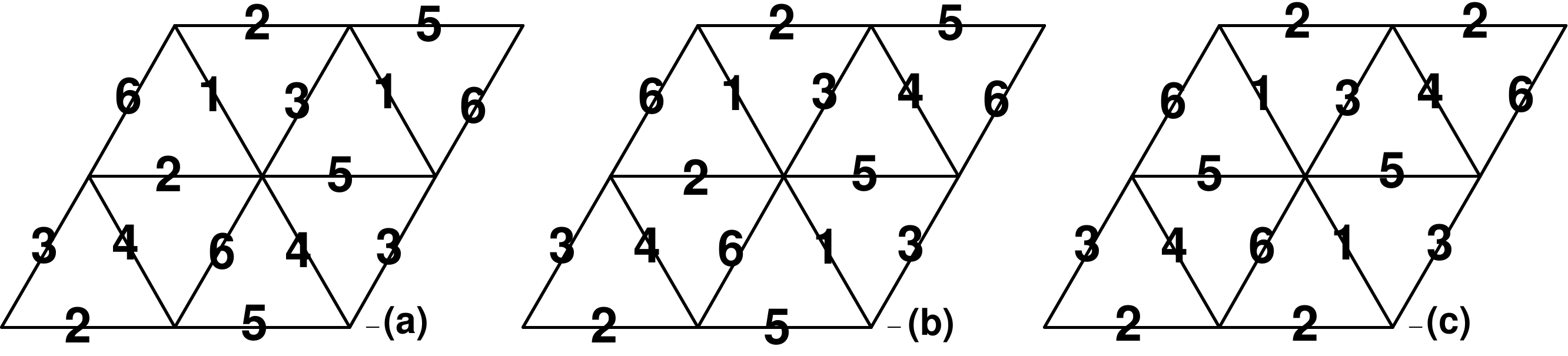}
  \vspace{-8pt}
  \caption{}
  \label{fig:23}
\end{figure}

\subsubsection{Figure 24 (a)}

In this figure we have a board of twelve equilateral parallelograms and the
plates are $\left[ 1316\right] $, $\left[ 2124\right] $, $\left[ 3235\right] 
$, $\left[ 4346\right] $, $\left[ 5154\right] $, $\left[ 6265\right] $,
twice.

\begin{figure}
  \centering
  \includegraphics[width=5.2719in]{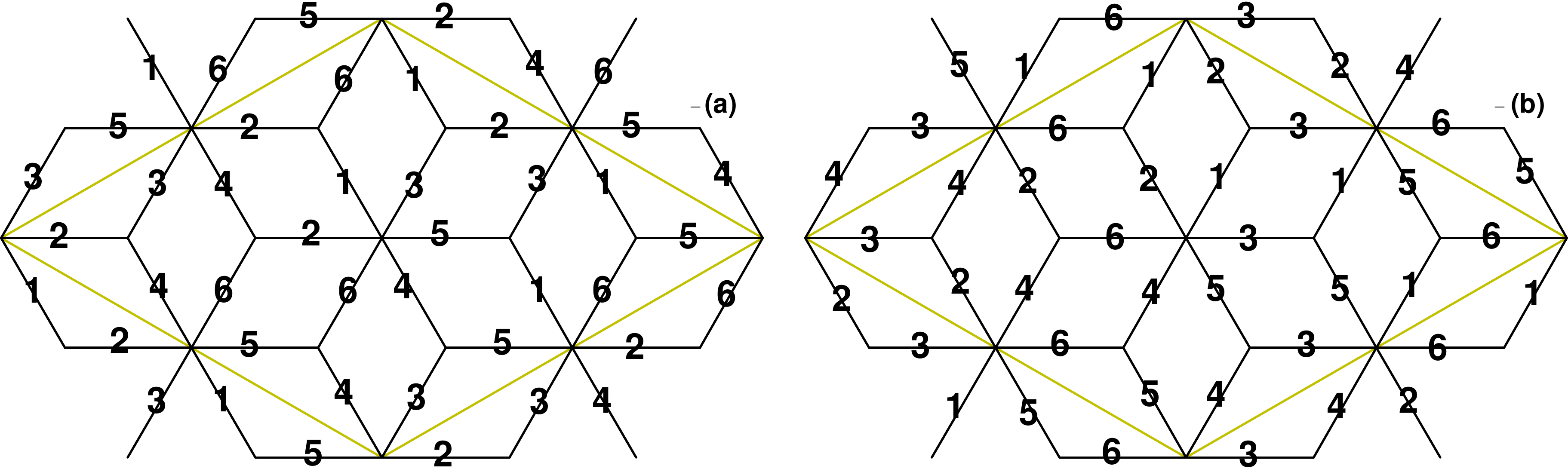}
  \caption{}
  \label{fig:24}
\end{figure}

\subsection{\boldmath Puzzles \guillemotleft $p6mm$\guillemotright}

\subsubsection{Figure 24 (b)}

The group is is generated by 
\begin{equation*}
\begin{tabular}{|c|c|c|c|}
\hline
$a$ & $d$ & $x$ & $y$ \\ \hline
$\left( 123456\right) $ & $\left( 156423\right) $ & $\left( 14\right) \left(
23\right) \left( 56\right) $ & $\left( 15\right) \left( 24\right) $ \\ \hline
\end{tabular}
\end{equation*}

In this figure we have a board of twelve equilateral parallelograms and the
plates are $\left[ 1313\right] $, $\left[ 1616\right] $, $\left[ 2121\right] 
$, $\left[ 2424\right] $, $\left[ 3232\right] $, $\left[ 3535\right] $,
$\left[ 4343\right] $, $\left[ 4646\right] $, $\left[ 5151\right] $, $\left[
5454\right] $, $\left[ 6262\right] $, $\left[ 6565\right] $.

\subsubsection{Figure 25}

In this figure, the board has $12$ faces ($4$ regular hexagons and $8$
equilateral triangles) and the plates are $\left[ 123456\right] $, $\left[
153426\right] $, $\left[ 156423\right] $, $\left[ 126453\right] $, $\left[
132\right] $, $\left[ 135\right] $, $\left[ 162\right] $, $\left[ 165\right] 
$, $\left[ 243\right] $, $\left[ 246\right] $, $\left[ 354\right] $, $\left[
465\right] $. This board and these plates are made in such a way that one
has a maximal solution with the group generated by 
\begin{equation*}
\begin{tabular}{|c|c|c|c|}
\hline
$a$ & $d$ & $x$ & $y$ \\ \hline
$\left( 123456\right) $ & $\left( 156423\right) $ & $\left( 14\right) \left(
23\right) \left( 56\right) $ & $\left( 15\right) \left( 24\right) $ \\ \hline
\end{tabular}
\end{equation*}

\begin{figure}
  \centering
  \includegraphics[width=3.8251in]{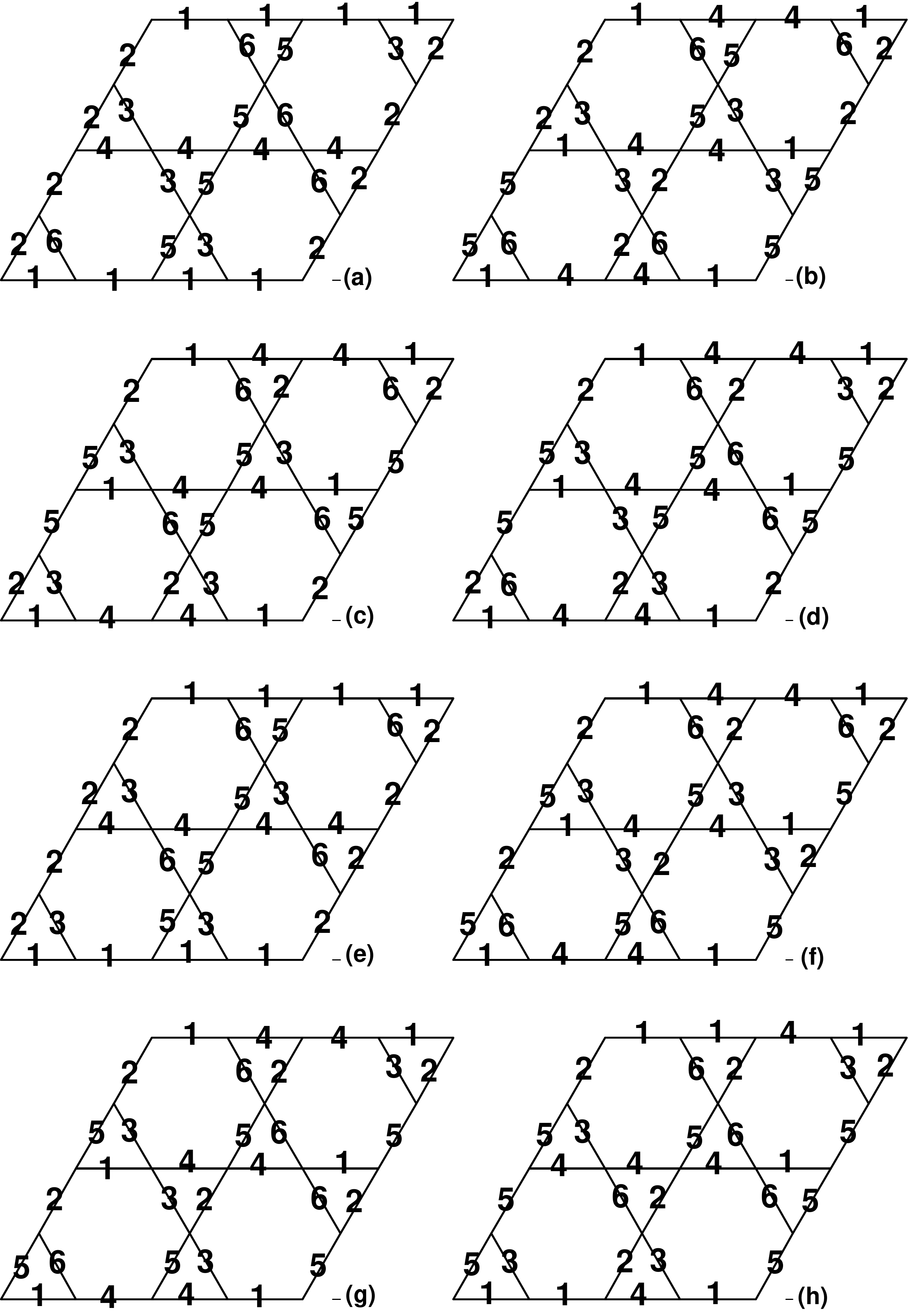}
  \vspace{-4pt}
  \caption{}
  \vspace{8pt}
  \label{fig:25}
\end{figure}

This puzzle puzzle has $8$ equivalence classes: (a) $1$ has a group of order 
$48$, (b) $1$ has a group of order $24$, (c)--(e) $3$ have a group of order
$16$, (f)--(g) $2$ have a group of order $8$, (e) $1$ has group of order $6$.
As $\left| G_P\right| =48$, one has that once we put a plate over a face
there are $32$ different possibilities ($32$ natural solutions) represented
in Figure 25: 
\begin{equation*}
48\left( \frac 1{48}+\frac 1{24}+\frac 3{16}+\frac 28+\frac 16\right) =32\text{.}
\end{equation*}

\begin{figure}
  \centering
  \includegraphics[width=2.1802in]{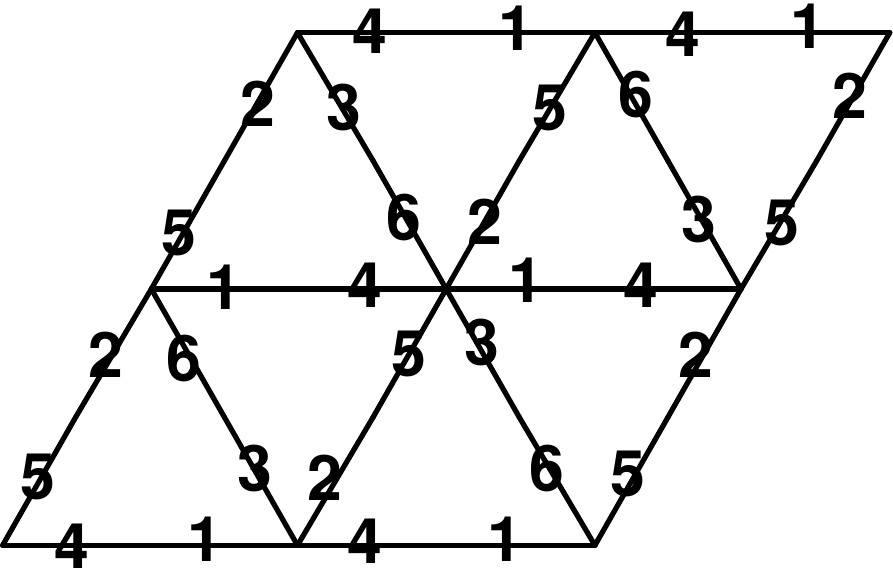}
  \vspace{-6pt}
  \caption{}
  \label{fig:26}
\end{figure}

\subsubsection{Figure 26}

The group is generated by 
\begin{equation*}
\begin{tabular}{|c|c|c|c|}
\hline
$a$ & $d$ & $x$ & $y$ \\ \hline
$\left( 123456\right) $ & $\left( 156423\right) $ & $\left( 26\right) \left(
35\right) $ & $\left( 12\right) \left( 36\right) \left( 45\right) $ \\ \hline
\end{tabular}
\end{equation*}

In this figure, the board has eight faces (equilateral triangles) and the
plates are $\left[ 14,36,25\right] $, $\left[ 14,36,52\right] $, $\left[
14,63,25\right] $, $\left[ 14,63,52\right] $, $\left[ 41,36,25\right] $,
$\left[ 41,36,52\right] $, $\left[ 41,63,25\right] $, $\left[ 41,63,52\right] 
$. Figure 26 shows a maximal solution that has all the $48$ symmetries of
the board.

\subsubsection{Figure 27 (a) and (b)}

In Figure 27 (a) we have a board of twelve equilateral parallelograms and
the plates are $\left[ 1245\right] $, $\left[ 1625\right] $, $\left[
1465\right] $, $\left[ 2136\right] $, $\left[ 2356\right] $, $\left[
3241\right] $, $\left[ 3461\right] $, $\left[ 4352\right] $, $\left[
5463\right] $.

In Figure 27 (b) we have a board of twelve equilateral parallelograms and
the plates are $\left[ 1221\right] $, $\left[ 1441\right] $, $\left[
1661\right] $, $\left[ 2332\right] $, $\left[ 2552\right] $, $\left[
3443\right] $, $\left[ 3663\right] $, $\left[ 4554\right] $, $\left[
5665\right] $.

The group of these maximal solutions of both puzzles is generated by 
\begin{equation*}
\begin{tabular}{|c|c|c|c|}
\hline
$a$ & $d$ & $x$ & $y$ \\ \hline
$\left( 123456\right) $ & $\left( 163254\right) $ & $\left( 14\right) \left(
23\right) \left( 56\right) $ & $\left( 15\right) \left( 24\right) $ \\ \hline
\end{tabular}
\end{equation*}

\begin{figure}
  \centering
  \includegraphics[width=4.3059in]{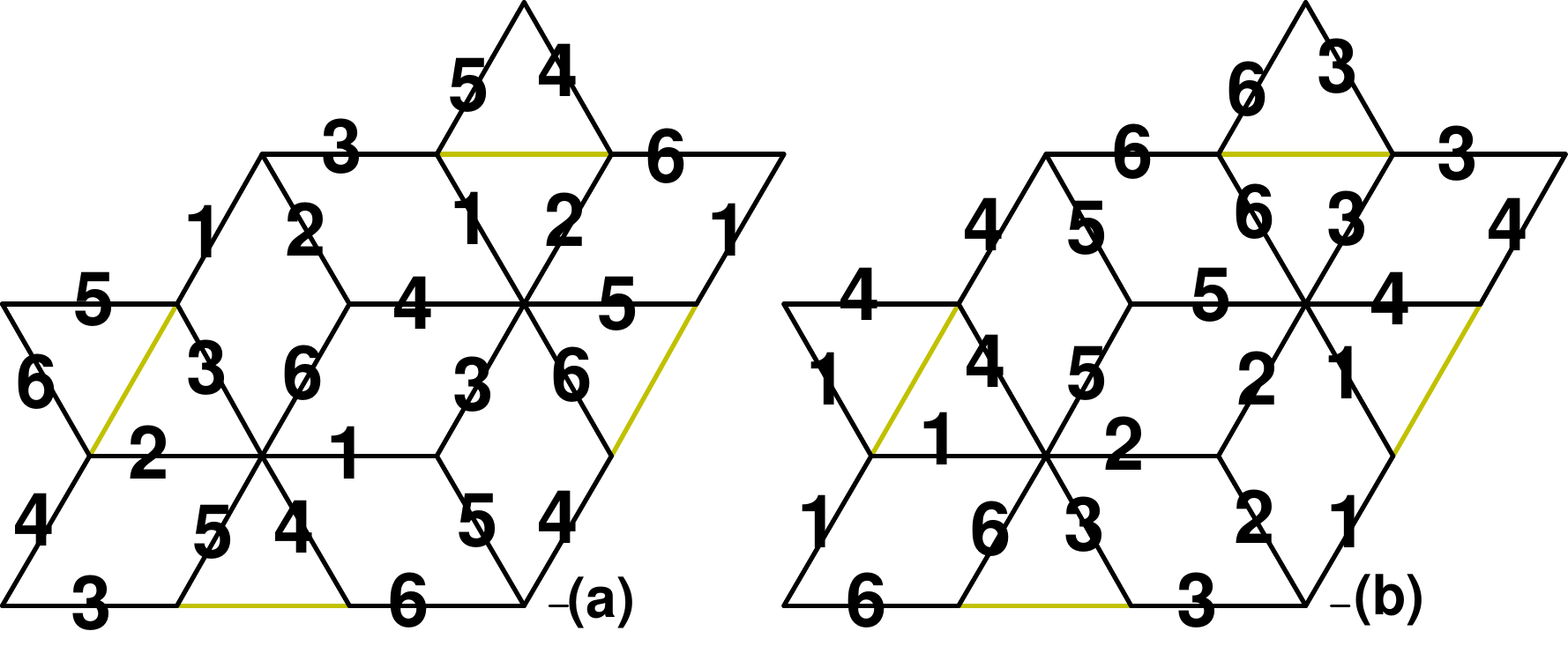}
  \vspace{-8pt}
  \caption{}
  \label{fig:27}
\end{figure}

\subsubsection{Figure 28}

The group is generated by 
\begin{equation*}
\begin{tabular}{|c|c|c|c|}
\hline
$a$ & $d$ & $x$ & $y$ \\ \hline
$\left( 123456\right) $ & $\left( 163254\right) $ & $\left( 14\right) \left(
23\right) \left( 56\right) $ & $\left( 15\right) \left( 24\right) $ \\ \hline
\end{tabular}
\end{equation*}

In this figure, the board has $9$ faces ($3$ regular hexagons and $6$
equilateral triangles) and the plates are: in (a), $\left[ 153\right] $,
$\left[ 264\right] $, three times, and $\left[ 123456\right] $, $\left[
163254\right] $, $\left[ 143652\right] $; in (b), $\left[ 111\right] $,
$\left[ 222\right] $ $\left[ 333\right] $, $\left[ 444\right] $ $\left[
555\right] $, $\left[ 666\right] $, $\left[ 123456\right] $, $\left[
163254\right] $, $\left[ 143652\right] $.

\begin{figure}
  \centering
  \includegraphics[width=5.2416in]{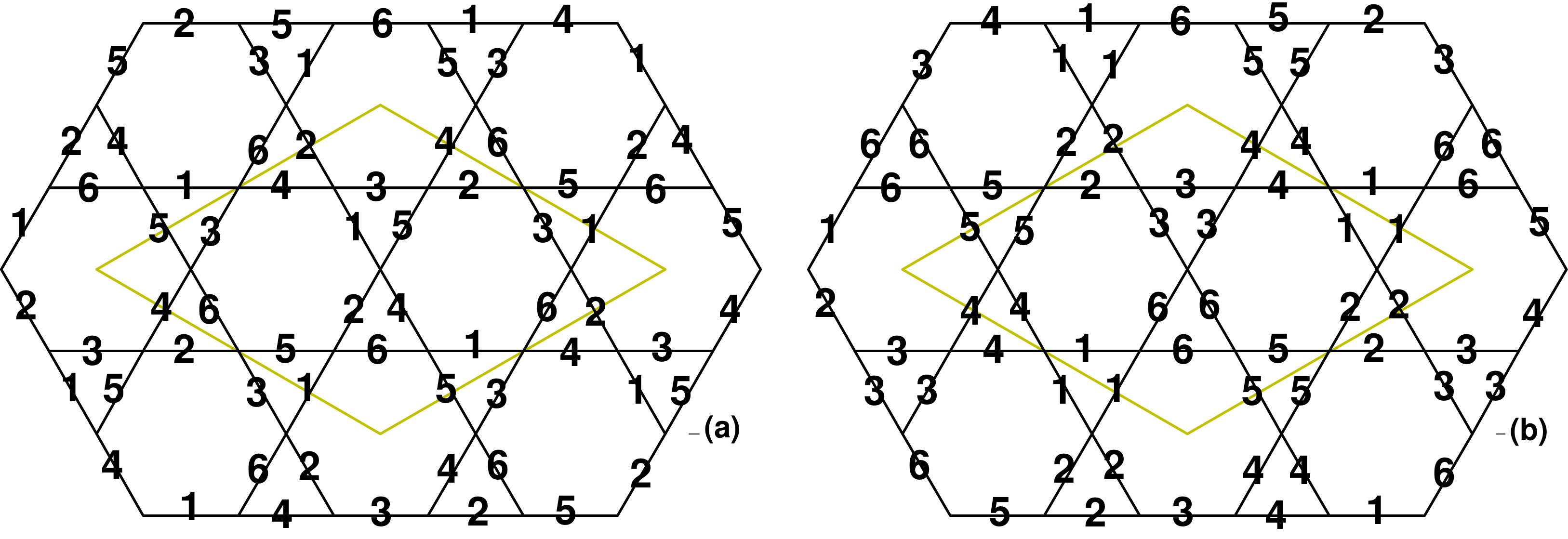}
  \vspace{-8pt}
  \caption{}
  \label{fig:28}
\end{figure}

\subsubsection{Figure 29 (a)}

The group is generated by 
\begin{equation*}
\begin{tabular}{|c|c|c|c|}
\hline
$a$ & $d$ & $x$ & $y$ \\ \hline
$\left( 123456\right) $ & $\left( 163254\right) $ & $\left( 14\right) \left(
23\right) \left( 56\right) $ & $\left( 15\right) \left( 24\right) $ \\ \hline
\end{tabular}
\end{equation*}

In this figure we have a board with $24$ triangular equilateral faces.

The plates are: $\left[ 153\right] $, $\left[ 264\right] $, three times;
$\left[ 124\right] $, $\left[ 125\right] $, $\left[ 132\right] $, $\left[
134\right] $, $\left[ 136\right] $, $\left[ 145\right] $, $\left[ 146\right] 
$, $\left[ 162\right] $, $\left[ 165\right] $, $\left[ 235\right] $, $\left[
236\right] $, $\left[ 243\right] $, $\left[ 245\right] $, $\left[ 256\right] 
$, $\left[ 346\right] $, $\left[ 354\right] $, $\left[ 356\right] $, $\left[
465\right] $.

\begin{figure}
  \centering
  \includegraphics[width=5.78in]{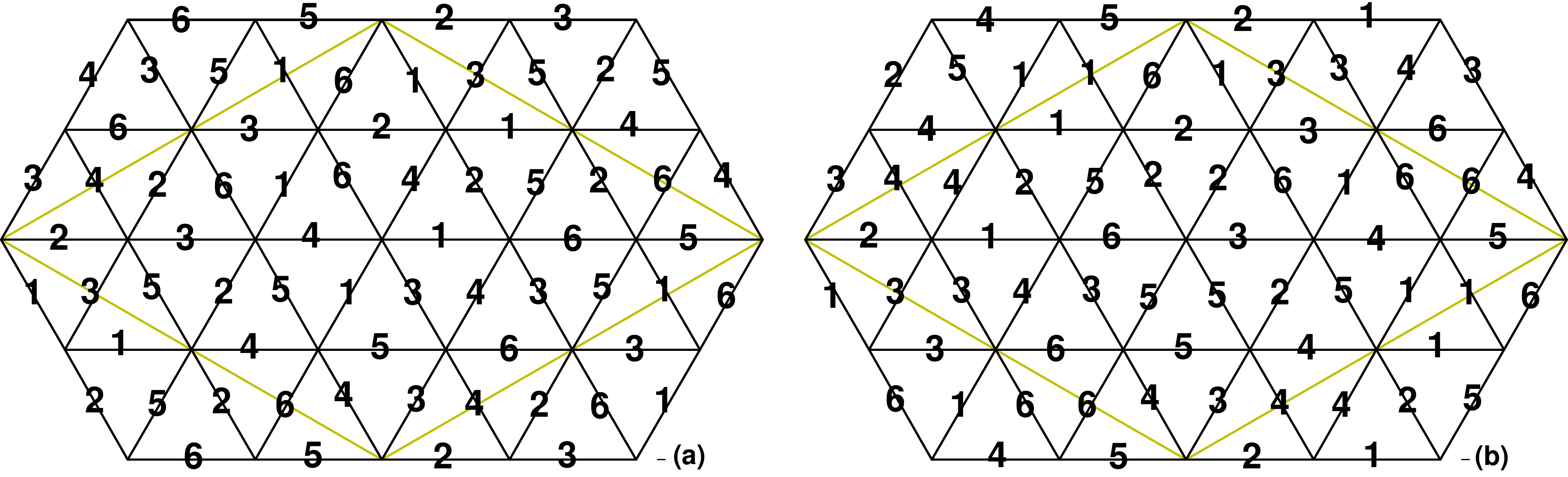}
  \vspace{-12pt}
  \caption{}
  \label{fig:29}
\end{figure}
\subsubsection{Figure 29 (b)}

The group is generated by 
\begin{equation*}
\begin{tabular}{|c|c|c|c|}
\hline
$a$ & $d$ & $x$ & $y$ \\ \hline
$\left( 123456\right) $ & $\left( 163254\right) $ & $\left( 14\right) \left(
23\right) \left( 56\right) $ & $\left( 15\right) \left( 24\right) $ \\ \hline
\end{tabular}
\end{equation*}

In this figure we have again a board with $24$ triangular equilateral faces.

The plates are: $\left[ 111\right] $, $\left[ 222\right] $ $\left[
333\right] $, $\left[ 444\right] $ $\left[ 555\right] $, $\left[ 666\right]$;
$\left[ 124\right] $, $\left[ 125\right] $, $\left[ 132\right] $, $\left[
134\right] $, $\left[ 136\right] $, $\left[ 145\right] $, $\left[ 146\right] 
$, $\left[ 162\right] $, $\left[ 165\right] $, $\left[ 235\right] $, $\left[
236\right] $, $\left[ 243\right] $, $\left[ 245\right] $, $\left[ 256\right] 
$, $\left[ 346\right] $, $\left[ 354\right] $, $\left[ 356\right] $, $\left[
465\right] $.


\bigskip

\begin{thebibliography}{9}
\bibitem{Armstrong}  M. A. Armstrong, \textit{Groups and Symmetry}. Berlin:
Springer-Verlag 1988.

\bibitem{Levine}  M. Levine, \textit{Plane symmetry groups}.\\
\url{http://www.math.uchicago.edu/~may/VIGRE/VIGRE2008/REUPapers/Levine.pdf}

\bibitem{rezende1}  J. Rezende, \textit{Puzzles with polyhedra and
permutation groups}.\\
\url{http://gfm.cii.fc.ul.pt/people/jrezende/jr_puzzles-poly-perm.pdf}

\bibitem{rezende2}  J. Rezende, \textit{On the Puzzles with polyhedra and
numbers} (2001).\\
\url{http://gfm.cii.fc.ul.pt/people/jrezende/jr_poliedros-puzzles_en.pdf}

\bibitem{rezende3}  J. Rezende, \textit{Puzzles with polyhedra and numbers},
(BGS Colloquium XI, April 2008), Actas, 103--131. Lisbon: LUDUS, 2009.\\
\url{http://ludicum.org/bgs08/bgs08-proceedings.pdf}

\bibitem{Schattschneider}  D. Schattschneider, \textit{The Plane Symmetry
Groups: Their Recognition and Notation}, American Mathematical Monthly,
Volume 85, Issue 6, 439--450.\\
\url{http://www.math.fsu.edu/~quine/MB_10/schattschneider.pdf}

\bibitem{Schwarzenberger}  R. L, Schwarzenberger, \textit{Colour symmetry}.
Bull.\ Lond.\ Math.\ Soc.\ 16, 209--240 (1984).

\bibitem{Wilson}  R. Wilson, P. Walsh, J. Tripp, I. Suleiman, R. Parker, S.
Norton, S. Nickerson, S. Linton, J. Bray, and R. Abbott, \textit{ATLAS of
Finite Group Representations --- Version 3.}\\
\url{http://brauer.maths.qmul.ac.uk/Atlas/v3/}
\end{thebibliography}
\end{document}